\documentclass[12pt]{amsart}
\usepackage{fullpage}
\usepackage{amscd,amsmath,amsthm,amssymb,graphics,color}

\usepackage[colorlinks,linkcolor=blue,filecolor=blue,citecolor=blue,urlcolor=blue,bookmarks=false,hypertexnames]{hyperref} 
\usepackage{tikz}
\usepackage{pdfpages}
\usepackage{graphicx}
\usepackage{paralist}
\usepackage{tkz-graph}
\usepackage{tikz-cd}
\usepackage{dirtytalk}
\usepackage{csquotes}

\usepackage{caption,subcaption}
\usepackage[capitalize]{cleveref}
\newcommand{\qand}{\quad \mbox{and} \quad}
\newcommand{\qor}{\quad \mbox{or} \quad}

\newcommand{\qforsome}{\quad \mbox{for some} \quad}
\newcommand{\qforall}{\quad \mbox{for all} \quad}
\newcommand{\qwith}{\quad \mbox{with} \quad}
\newcommand{\qwhere}{\quad \mbox{where} \quad}
\newcommand\reals{\mathbb{R}}
\newcommand{\MA}{\mathcal A}

\def\MH{{\mathcal H}}

\def\MP{{\mathcal P}}
\def\MQ{{\mathcal Q}}

\def\opn#1#2{\def#1{\operatorname{#2}}} % to make operators
\opn\det{det}
\opn\ini{in}
\opn\height{height}
\opn\rank{rank}
\opn\supp{supp}
\opn\vsupp{vsupp}
\opn\gcd{gcd}
\opn\Min{Min}
\opn\sp{sp}
\opn\lk{lk}
\opn\st{st}
\opn\LCM{LCM}
\opn\lcm{lcm}

\DeclareMathOperator{\im}{im}
\DeclareMathOperator{\Taylor}{Taylor}
\DeclareMathOperator{\Scarf}{Scarf}

\makeatletter
\@addtoreset{equation}{section}
\def\theequation{\thesection.\@arabic \c@equation}
\makeatother

\theoremstyle{plain}

\newtheorem{theorem}[equation]{Theorem}

\newtheorem{lemma}[equation]{Lemma}
\newtheorem{proposition}[equation]{Proposition}

\theoremstyle{definition}

\newtheorem{definition}[equation]{Definition}

\newtheorem{discussion}[equation]{Discussion}

\newtheorem{notation}[equation]{Notation}

\newtheorem{Remark}[equation]{Remark}
\newenvironment{remarkbox}[1][]{%
\begin{Remark}[#1]\pushQED{\qed}}{\popQED \end{Remark}}
\newtheorem{example}[equation]{Example}

\newtheorem{setup}[equation]{Setup}
\address{Department of Mathematics \& Statistics,
 Dalhousie University,
 6297 Castine Way,
 PO BOX 15000,
 Halifax, NS,
 Canada B3H 4R2
 } 
\email{
l.bu@dal.ca, 
faridi@dal.ca, 
md472930@dal.ca, 
hollebenthiago@dal.ca,
mh290120@dal.ca,
d.veer@dal.ca,
yp259112@dal.ca,
Scott.Wesley@dal.ca
}
\author{Louis Bu} 
\author[S. Faridi]{Sara Faridi}
\author{Iresha Madduwe Hewalage} 
\author{Thiago Holleben} 
\author{Hasan Mahmood} 
\author{Dharm Veer} 
 \author{Kyle Wang} 
\author{Scott Wesley} 

\title{When are Morse resolutions polyhedral?}
\begin{document}
\maketitle

\begin{abstract}
 It is known that the chain complex of a simplex on $q$ vertices   can be used to construct a free resolution of any ideal generated by $q$ monomials, and as a direct result, the Betti numbers always have binomial upper bounds, given by the number of faces of a simplex in each dimension.
 
 It is also known that for most monomials the resolution provided by the simplex is far from minimal. Discrete Morse theory provides an algorithm  called \say{Morse matchings} by which faces of the simplex can be removed so that the chain complex on the remaining faces is still a free resolution of the same ideal. An immediate positive effect is an often considerable improvement on the bounds on Betti numbers.  A caveat is the loss of the combinatorial structure of the simplex we started with:  the output of the Morse matching process is a cell complex with no obvious structure besides an \say{address} for each cell. 

The main question in this paper is: which Morse matchings lead to Morse complexes that are polyhedral cell complexes?  We prove that if a monomial ideal is minimally generated by up to four generators, then there is a maximal Morse matching of the simplex such that the resulting cell complex is a polyhedral cell complex. We then give an example of a monomial ideal minimally generated by six generators whose minimal free resolution is supported on a Morse complex and the Morse complex cannot be polyhedral no matter what Morse matching is chosen, and we go further to show that this ideal cannot have any polyhedral minimal free resolution.
\end{abstract}

%%%%%%%%%%%%%%%%%%%%%%%%%%%%%%%%%%%%%%%%%%%%%%%%%%%%%%%%%%%%%%%%%%%%%%%%%
\section{Introduction}
%%%%%%%%%%%%%%%%%%%%%%%%%%%%%%%%%%%%%%%%%%%%%%%%%%%%%%%%%%%%%%%%%%%%%%%%%

Discrete Morse theory provides a powerful tool to analyze free resolutions of monomial ideals via geometric objects. A free resolution is a way to encode relations between homogeneous polynomials with coefficients in a field into an exact sequence of free module homomorphisms.  The smallest such sequence, in terms of the number and ranks of the free modules in the sequence, is called a {\it minimal free resolution}, and is  unique up to isomorphism of complexes. 

In 1966, Taylor~\cite{Taylorthesis} began the study of resolutions from a combinatorial perspective, by using algebraic combinatorics to compute a free resolution of  monomial ideals. She showed that the simplicial chain complex of a simplex with $r$ vertices can be \enquote{homogenized} to give a free resolution of a monomial ideal with $r$ minimal generators. This free resolution is known as the \emph{Taylor resolution}, which is most often  non-minimal. 

The term \enquote{homogenization} was first used by Bayer and Sturmfels~\cite{BS98cellularresolutions}, who generalized Taylor's construction to compute free resolutions of monomial ideals by homogenizing the chain complex of cell complexes.
 When the homogenization of the chain complex of a cell complex $\Gamma$ gives rise to a free  resolution of a monomial ideal $I$, we say $\Gamma$ \emph{supports} a free resolution  of $I$, and this free resolution is called a \emph{cellular resolution}.

Every monomial ideal has a cellular resolution -- the Taylor resolution -- but 
Reiner and Welker~\cite{Reiner_Welker_2001} and Velasco~\cite{Velasco08minimalfreeresolution} proved that there are classes of monomial ideals whose minimal free resolution cannot be supported by \emph{any} simplicial, respectively even cellular, complex. Nevertheless, geometrically structured cellular resolutions such as simplicial or polyhedral ones which are (close to) 
minimal give us invaluable information about the free resolution of the ideal.

To compute cellular resolutions which are closer to minimal ones than the  Taylor resolution, a natural method is to adapt tools from  discrete homotopy theory   to shrink the size of the cell complex supporting a resolution of an ideal.
Discrete Morse theory~(\cite{BW2002cellularresolution}) is one such tool and the premise of our paper. The idea is to begin with the face poset of the Taylor complex and systematically eliminate pairs of faces using \say{Morse matchings}.

The output of a Morse matching is a \say{Morse complex}, whose $i$-faces are in one to one correspondence with the unmatched faces of the Taylor complex. While the main theorem in this area states that the resulting Morse complex supports a resolution of the ideal -- giving an excellent bound on the size of the minimal free resolution -- the structure of the Morse complex itself can  be  elusive.

In this article, we give an example of a square-free monomial ideal $I$ with six minimal generators, whose minimal free resolution is supported on a Morse complex, and the Morse complex cannot be polyhedral no matter what Morse matching is chosen. 

The \emph{Scarf complex} of a monomial ideal is a simplicial complex whose homogenized chain complex appears inside any free resolution of that ideal. 
Notably for our example, the Morse complex differs from the Scarf complex by only one $3$-face. The Morse complex being non-polyhedral in this case makes that one extra face attach to the Scarf complex \say{along} non-faces.

%We characterize all maximal Morse matchings of the face poset of the Taylor complex of $I$. For each maximal Morse matching, the resulting cell complex (in the sense of~\cite{BW2002cellularresolution}) is non-polyhedral and supports the minimal free resolution of $I$.

We also show that when a monomial ideal $I$ is minimally generated by four elements, that is, when the Taylor complex has dimension $\leq 3$, there is always a maximal Morse matching of the Taylor complex of $I$ such that the resulting cell complex is a polyhedral cell complex. We do not know if the polyhedral Morse resolution for these ideals is minimal, but it is already known that monomial ideals with up to five generators  have minimal cellular resolutions~\cite{Roya,chau2025minimal,montaner2025morse}. 

This paper is organized as follows: In \cref{sec:polyhedral}, we define polyhedral cell complexes, review multigraded Betti numbers, Taylor and Scarf complexes, and cellular resolutions of monomial ideals. 
In \cref{sec:morsetheory}, we recall the definitions from discrete Morse theory and results regarding its application in computing cellular resolutions of monomial ideals. \cref{sec:4generators} is dedicated to showing that for a monomial ideal minimally generated by four elements, there exists a maximal homogeneous acyclic matching such that the resulting cell complex is a polyhedral cell complex. In \cref{sec:non-polyhedralexample}, we then give an example of a monomial ideal with $6$ minimal generators whose minimal Morse resolution is not supported on a polyhedral cell complex. We conclude by showing in \cref{r:minimal} that this ideal cannot have any minimal polyhedral free resolution.  Our example is geometrically optimal, as only one multidegree in the multigraded free resolution of this ideal is not a Scarf multidegree, and therefore the problem reduces to one of polyhedral geometry: there is no way to attach a simplicial $3$-polytope with the correct multidegree to a certain simplicial complex so that the resulting complex is a polyhedral cell complex. 
We would like to point out that an example of a monomial ideal with $13$ generators and no polyhedral minimal resolution was given in~\cite{Reiner_Welker_2001}.

\subsection*{Acknowledgements}
We are grateful to Linh Dinh for her contributions at the earlier stages of this project, and to Volkmar Welker for helpful comments. Faridi's research is supported by  NSERC Discovery Grant~2023-05929. Veer is supported by an AARMS Postdoctoral Fellowship. Wesley is supported by a Level~2 Killam Predoctoral Scholarship. 

%%%%%%%%%%%%%%%%%%%%%%%%%%%%%%%%%%%%%%%%%%%%%%%%%%%%%%%%%%%%%%%%%%%%%%%%%
\section{Simplicial, cellular, and polyhedral resolutions}\label{sec:polyhedral}
%%%%%%%%%%%%%%%%%%%%%%%%%%%%%%%%%%%%%%%%%%%%%%%%%%%%%%%%%%%%%%%%%%%%%%%%%

Throughout this paper, we will deal with three different notions of combinatorial objects: simplicial, polyhedral and cell complexes. We begin by recalling their definitions, and we note that we often simplify definitions to avoid introducing extra  terminology.

\subsection{Cell complexes}

\begin{definition}[{\bf convex polytopes, simplices}]\label{d:polytope} A \emph{convex polytope} (see~\cite{GS69convexpolytopes}) $\MP$ in  $\reals^N$ is the convex hull of a finite set of points $\{a_0,\ldots, a_n\}$ in  $\reals^N$
\[
\MP = \big \{x\in \reals^N: x = \sum_{i=0}^{n} t_ia_i, \quad \text{where} \quad \sum_{i=0}^{n} t_i =1 \qand t_i\geq 0 \qforall i \big \}.
\]
The \emph{dimension} of $\MP$, denoted by $\dim(\MP)$ is  $d-1$ where $d$ is the vector space dimension of the space spanned by $\{a_0,\ldots, a_n\}$. 
A $d$-polytope is a convex polytope of dimension $d$.

If  the vectors $a_1-a_0, \ldots, a_n-a_0$ are linearly independent in $\reals^N$,  then $\MP$ is an $n$-dimensional convex polytope called  an \emph{$n$-simplex}   (see~\cite{munkresalgtop}). 
When $N=n+1$ and $\{a_0, \ldots, a_n\}$ is the standard unit basis of $\reals^{n+1}$, 
the $n$-simplex spanned by $a_0,\ldots,a_n$ is called the \emph{standard $n$-simplex}, denoted by $\Delta_n$.   
\end{definition}

A hyperplane $\MH$ \emph{supports} a convex polytope $\MP$ if $\MP\cap \MH \neq \emptyset$ and $\MP$ lies in one of the two closed half-spaces bounded by $\MH$. If $\MH$ supports $\MP$, then $\MP\cap \MH$ is called a \emph{face} of $\MP$. For convenience, to a $d$-polytope $\MP$ we add two faces: $\emptyset$ of dimension $-1$ and $\MP$ itself, of dimension $d$.

For a convex polytope $\MP$, the following properties hold:
\begin{enumerate}
    \item 
        Every face of $\MP$ is a polytope.
    \item     
        Every face of a face of $\MP$ is also a face of $\MP$.
    \item 
        For every two faces $F, G$ of $\MP$, $F\cap G$ is also a face of $\MP$. We write $F\wedge G$ for $F\cap G$.
    \item 
        For every two faces $F, G$ of $\MP$, there exists a unique \say{smallest} face $F\vee G$ of $\MP$, which contains both $F$ and $G$.
\end{enumerate}

The set of all faces of $\MP$ forms a lattice with respect to the operations $\wedge$ and $\vee$ defined above, called the \emph{face-lattice} of $\MP$.

A polyhedral complex is a collection of convex polytopes glued along faces ({\cite[Definition~4.1]{millersturmfels}}, {\cite[\S2]{munkresalgtop}})

\begin{definition}[{\bf polyhedral and simplicial complexes}]
A \emph{polyhedral cell complex} $X$ is a collection of convex polytopes, called \emph{faces} of $X$, such that
\begin{enumerate}
    \item every face of a polytope in $X$ is also a face of $X$, and
    
    \item  the intersection of any two faces $\MP$ and $\MQ$ of $X$ is a face of both $\MP$ and $\MQ$.
\end{enumerate}

When all the convex polytopes in a polyhedral cell complex $X$ are simplices, $X$ is called a \emph{simplicial complex}.
\end{definition}

Polyhedral cell complexes can be defined in an abstract way, independent of their embedding into Euclidean space. In the case of simplicial complexes, the abstract definition is the one that is most commonly used, as every abstract simplicial complex is geometrically realizable. For polytopes on the other hand, the realizability problem of an abstract polytope as a convex polytope is very hard, and is often called \emph{Steinitz problem} (see for example~\cite{BS1989}), due to Steinitz' celebrated theorem~\cite{zieglerbook}.

\begin{definition}[{\bf abstract simplicial complex} ]
An \emph{abstract simplicial complex} $\Delta$ on vertex set $[n] = \{1,\dots, n\}$ is a collection of subsets of $[n]$ called \emph{faces} such that $\sigma \in \Delta$ and $\tau \subset \sigma$ implies $\tau \in \Delta$.
The \emph{dimension} of a face $\sigma$ in $\Delta$, denoted by $\dim(\sigma)$, is the number $|\sigma|-1$. The \emph{dimension} of the simplicial complex $\Delta$, denoted by $\dim(\Delta)$, is $\max \{\dim(\sigma) : \sigma \in \Delta\}$. The maximal faces under inclusion are called \emph{facets} of $\Delta$. An $n$-dimensional simplicial complex with exactly one facet is called a \emph{$n$-simplex}. 
\end{definition}

Denote by $\partial\Delta_n$ the \emph{boundary} of standard $n$-simplex $\Delta_n$, (\cref{d:polytope}) which is the union of all $(n-1)$-dimensional faces of $\Delta_n$.
An $n$-dimensional \emph{cell complex} $X$ is the last space in an inductively constructed sequence of topological spaces $X^{(0)} \subseteq  X^{(1)} \subseteq \cdots \subseteq X^{(n)}=X$, following the procedure outlined below:
\begin{enumerate}
    \item Start with a discrete set $X^{(0)}$, whose points are regarded as $0$-cells.
    \item For each integer $0 <k \leq n$, we inductively construct the $k$-skeleton $X^{(k)}$ from $X^{(k-1)}$ by attaching $k$-simplices via gluing maps $\varphi : \partial \Delta_k \to X^{(k-1)}$ that maps the interior of each face of $\Delta_k$ homeomorphically to the interior of a simplex in $X^{(k-1)}$ of the same dimension. Thus, $\varphi$ maps each vertex of $\Delta_k$ to a point in $X^{(0)}$, each edge of $\Delta_k$ to an edge in $X^{(1)}$, and so on; More formally, $X^{(k)}$ is defined as the quotient space
    \[
    X^{(k)} := \Big(X^{(k-1)} \sqcup \bigsqcup_{i\in I_k} \Delta_{k,i} \Big) / \sim,
    \]
    for some index set $I_k$, where $\sqcup$ denotes disjoint union, each $\Delta_{k,i}$ is a disjoint copy of $\Delta_k$, and $x \sim \varphi_i(x)$ for each gluing map $\varphi_i : \partial \Delta_{k,i} \to X^{(k-1)}$. The interiors of the simplices $\Delta_{k,i}$ are the $k$-cells of $X$.
\end{enumerate}

Every simplicial complex is a polyhedral cell complex; however, the converse does not hold. For example, the convex polytope shown below on the left is a polyhedral cell complex, but it is not a simplicial complex as $\{1,2,3,4\}$ is not a simplex.
Every simplicial and polyhedral cell complex is indeed a cell complex (see \cite[Section~2.2]{kozlovalgtop}) but, again,  the converse is not true. For example, the cell complex below on the right is not a polyhedral cell complex.
   \[
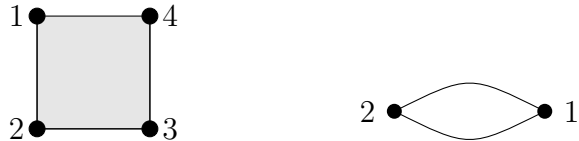
   \begin{tabular}{cc}
         \begin{tikzpicture}[scale = 1.5]
   \tikzstyle{point}=[inner sep=0pt]
   \node (a)[point,label=left:{$1$}] at (0,1) {}; 
   \node (b)[point,label=left:{$2$}] at (0,0) {}; 
   \node (c)[point,label=right:{$3$}] at (1,0) {}; 
   \node (d)[point,label=right:{$4$}] at (1,1) {}; 
   \draw [fill=gray!20](a.center) -- (b.center) -- (c.center) -- (d.center);
   \draw (a.center) -- (b.center);
   \draw (b.center) -- (c.center);
   \draw (a.center) -- (d.center);
    \draw (c.center) -- (d.center);
\foreach \p in {a,b,c,d}  \filldraw[black] (\p) circle (2pt);
         \end{tikzpicture}
   \qquad & \qquad \qquad
         \begin{tikzpicture}[scale = 1]
   \tikzstyle{point}=[inner sep=0pt]
   \node (a)[point, circle, fill=black, inner sep=2pt, label=right:{$1$}] at (2,1) {}; 
   \node (b)[point,circle, fill=black, inner sep=2pt, label=left:{$2$}] at (0,1) {}; 

   \draw (b.center) .. controls +(1,.5) .. (a.center);
   \draw (b.center) .. controls +(1,-.5) .. (a.center);
         \end{tikzpicture}
   \end{tabular}
   \]

In the following subsection, we utilize the previously discussed containment to differentiate the set of monomial ideals in a polynomial ring over a field. This distinction is based on whether the minimal free resolution of a monomial ideal is supported on a simplicial, polyhedral, or a cellular complex.

%simplicial $\subset$ polyhedral $\subset$ cellular complex
%%%%%%%%%%%%%%%%%%%%%%%%%%%%%%%%%%%%%%%%%%%%%%%%%%%%%%%%%%%%%%%%%%%%%%%%%
\subsection{Free resolutions}\label{sec:freeresolutions}
%%%%%%%%%%%%%%%%%%%%%%%%%%%%%%%%%%%%%%%%%%%%%%%%%%%%%%%%%%%%%%%%%%%%%%%%%

Let $S=K[x_1,\ldots,x_n]$ be a polynomial ring in $n$ variables over a field $K$, and $I=(m_1,\ldots,m_r)$ be a monomial ideal of $S$ minimally generated by $r$ monomials. 
The \emph{lcm lattice} of $I$, denoted by $\LCM(I)$, is the poset consisting of the least common multiples  of the generating set of $I$, ordered by divisibility. 

A \emph{graded free resolution} of $I$ over $S$ is an exact sequence of graded free $S$-modules of the form:
 \[
 \mathbb{F} : 0 \to F_d \stackrel{\partial_d}{\longrightarrow} \cdots \to F_i \stackrel{\partial_i}{\longrightarrow} F_{i-1} \to \cdots F_1 \stackrel{\partial_1}{\longrightarrow} F_0
 \]
where $I\cong F_0/\im(\partial_1)$, and each map $\partial_i$ is a homogeneous map of degree zero, i.e, a degree-preserving map.
We say that the free resolution $\mathbb{F}$ is \emph{minimal} if $\partial_i(F_i)\subseteq (x_1,\ldots,x_n)F_{i-1}.$
The minimal free resolution of $I$ is unique up to isomorphism of complexes. When $\mathbb{F}$ is the minimal free resolution of $I$, the rank of $F_i$ over $S$ is the $i$-th \emph{Betti number} of $I$, denoted by $\beta_i.$ 
The grading on each of the free modules $F_i$ can be refined by writing $F_i$ as direct sum of 1-dimensional free $S$-modules $S(m)$ indexed by monomials $m$ in the lcm lattice of $I$.
When $\mathbb{F}$ is a minimal resolution, one has
\[
F_i\cong \mathop{\bigoplus}_{m\in \LCM(I)} S(m)^{\beta_{i,m}}
\]
where the number $\beta_{i,m}$ is the $i$-th \emph{multigraded Betti number} of $I$ in multidegree $m.$

The \emph{Taylor complex} of $I$, denoted by $\Taylor(I)$ is a simplex on $r$ vertices, each vertex of the simplex is labeled with one of the monomials $m_1,\ldots,m_r$, and then each face is labeled by the lcm of the monomial labels of its vertices. 
%The monomial labels of of the Taylor complex belongs to the $\LCM(I)$. 
Taylor~\cite{Taylorthesis} proved that the simplicial chain complex of Taylor(I) %\dharm{can be \say{homogenized} to} 
gives rise to a multigraded free resolution of $I$.
This free resolution is known as the \emph{Taylor resolution}, which is typically non-minimal.
%Taylor’s resolution is usually non-minimal.
Every monomial ideal has a multigraded minimal free resolution contained in the Taylor resolution.

Following  the language of Bayer and Sturmfels~\cite{BS98cellularresolutions}, when the \say{homogenization} of a cell complex $\Gamma$ gives rise to a free  resolution of $I,$ we say $\Gamma$ \emph{supports a free resolution} of $I$.

A monomial ideal $I$ which has a (minimal) free resolution supported on a simplicial, polyhedral, or more generally cellular complex $\Gamma$ is said to have, respectively, a (minimal) {\it simplicial},  {\it polyhedral}, or {\it cellular}  free resolution. Note that both simplicial and polyhedral resolutions are also cellular resolutions. Moreover, we can check minimality of the resolution via the following criterion:
\begin{equation}\label{eq:minimal} 
\mbox{the resolution supported on $\Gamma$ is minimal} \iff  
m_\sigma \neq m_\tau \qforall \sigma, \tau \in \Gamma \qwith \sigma \prec \tau .
\end{equation}

Velasco~\cite{Velasco08minimalfreeresolution} proved that there are classes of monomial ideals whose minimal free resolution cannot be supported by any cell complex.
The \emph{Scarf complex} of I, denoted by $\Scarf(I)$, is a subcomplex of $\Taylor(I)$ obtained by removing
all faces of $\Taylor(I)$ that share a monomial label with another face. 
%The relevance of Scarf complexes comes from the fact that a minimal free resolution of $I$ contains an isomorphic copy of....
If a cell complex $\Gamma$ supports a free resolution of $I$, then $\Scarf(I)$ is a subcomplex of $\Gamma$~\cite[Proposition\ 59.4]{peevagradedsyzygies}. 
%$\Scarf(I)$ is a subcomplex of every cell complex that supports a free resolution of $I$.
When $\Scarf(I)$ is acyclic, it supports a free resolution of $I$ which is necessarily minimal.

\begin{example}\label{example}
 Let $I= (x,y,z)^2 =(x^2,xy,y^2,yz,z^2,xz) \subset S=K[x,y,z]$ be the monomial ideal.  The simplicial complex $\Scarf(I)$ is the simplicial complex on the right in \cref{fig:example}. Since $\Scarf(I)$ is not acyclic, it does not support a free resolution of $I$, but we know it shows up in any free resolution of $I$. It has been shown in~\cite{Faridi21secondpowerresolution} that $I$ has a free resolution supported on the  simplicial complex $\Delta$ which appears in \cref{fig:example} on the left. Using the minimality criterion in \eqref{eq:minimal}, we see that the central triangle and any one of each edges share the monomial label $xyz$, and so the resolution is not minimal.  Indeed a computation by Macaulay2~\cite{M2} shows that  $\beta_1(I) =8$ and $\beta_2(I)=3$, so one of the $9$ edges of $\Delta$ and one of the triangles is extra. Our argument above, and also the fact that all outside edges are Scarf, lead us to conclude  that the candidates for the extra edge are the three central  edges with monomial label $xyz$. 

   \begin{figure}[hbt!]
   \begin{tabular}{ccc}
         \begin{tikzpicture}[scale = 1]
   \tikzstyle{point}=[inner sep=0pt]
   \node (a)[point,label=above:{ $xy$}] at (0,1) {}; 
   \node (b)[point,label=left:{ $xz$}] at (-1,0) {}; 
   \node (c)[point,label=right:{ $yz$}] at (1,0) {}; 
   \node (d)[point,label=left:{ $x^2$}] at (-1.5,1.5) {}; 
   \node (e)[point,label=right:{ $y^2$}] at (1.5,1.5) {}; 
   \node (f)[point,label=right:{ $z^2$}] at (0,-1) {}; 
   \draw [fill=gray!20](a.center) -- (b.center) -- (c.center);
   \draw [fill=gray!20](a.center) -- (b.center) -- (d.center);
   \draw [fill=gray!20](a.center) -- (c.center) -- (e.center);
   \draw [fill=gray!20](b.center) -- (c.center) -- (f.center);
   \draw (a.center) -- (b.center);
   \draw (a.center) -- (c.center);
   \draw (a.center) -- (d.center);
   \draw (a.center) -- (e.center);
   \draw (b.center) -- (c.center);
   \draw (b.center) -- (d.center);
   \draw (b.center) -- (f.center);
   \draw (c.center) -- (e.center);
   \draw (c.center) -- (f.center);
   \node [label={ $xyz$}] at (0,0) {}; 
   \node [label={ $x^2yz$}] at (-.8,.4) {}; 
   \node [label={ $xy^2z$}] at (0.8,.4) {}; 
   \node [label={ $xyz^2$}] at (0,-.8) {}; 
  \foreach \p in {a,b,c,d,e,f}  \filldraw[black] (\p) circle (2pt);
  \end{tikzpicture}
   \qquad & \qquad 
         \begin{tikzpicture}[scale = 1]
   \tikzstyle{point}=[inner sep=0pt]
   \node (a)[point,label=above:{ $xy$}] at (0,1) {}; 
   \node (b)[point,label=left:{ $xz$}] at (-1,0) {}; 
   \node (c)[point,label=right:{ $yz$}] at (1,0) {}; 
   \node (d)[point,label=left:{ $x^2$}] at (-1.5,1.5) {}; 
   \node (e)[point,label=right:{ $y^2$}] at (1.5,1.5) {}; 
   \node (f)[point,label=right:{ $z^2$}] at (0,-1) {}; 
   %\draw [fill=gray!20](a.center) -- (c.center) -- (b.center)--(d:center);
   \shade (a.center) -- (b.center) -- (d.center);
   \draw [fill=gray!20](a.center) -- (c.center) -- (e.center);
   \draw [fill=gray!20](b.center) -- (c.center) -- (f.center);
   %\draw (a.center) -- (b.center);
   \draw (a.center) -- (c.center);
   \draw (a.center) -- (d.center);
   \draw (a.center) -- (e.center);
   \draw (b.center) -- (c.center);
   \draw (b.center) -- (d.center);
   \draw (b.center) -- (f.center);
   \draw (c.center) -- (e.center);
   \draw (c.center) -- (f.center);
   %\node [label={ $xyz$}] at (0,0) {}; 
   \node [label={ $x^2yz$}] at (-.8,.4) {}; 
   \node [label={ $xy^2z$}] at (0.8,.4) {}; 
   \node [label={ $xyz^2$}] at (0,-.8) {}; 
  \foreach \p in {a,b,c,d,e,f}  \filldraw[black] (\p) circle (2pt);
  \end{tikzpicture}
   \qquad & \qquad 
        \begin{tikzpicture}[scale = 1]
   \tikzstyle{point}=[inner sep=0pt]
   \node (a)[point,label=above:{ $xy$}] at (0,1) {}; 
   \node (b)[point,label=left:{ $xz$}] at (-1,0) {}; 
   \node (c)[point,label=right:{ $yz$}] at (1,0) {}; 
   \node (d)[point,label=left:{ $x^2$}] at (-1.5,1.5) {}; 
   \node (e)[point,label=right:{ $y^2$}] at (1.5,1.5) {}; 
   \node (f)[point,label=right:{ $z^2$}] at (0,-1) {}; 
   \draw (a.center) -- (d.center);
   \draw (a.center) -- (e.center);
   \draw (b.center) -- (d.center);
   \draw (b.center) -- (f.center);
   \draw (c.center) -- (e.center);
   \draw (c.center) -- (f.center);
  \foreach \p in {a,b,c,d,e,f}  \filldraw[black] (\p) circle (2pt);
  \end{tikzpicture}\\
         $\Delta$ && $\Scarf(I)$
          \end{tabular}
     \caption{Figure for \cref{example}}\label{fig:example}
   \end{figure}

If we remove the central triangle and one of the edges, we will have the middle picture in \cref{fig:example}, whose number of faces match  the Betti numbers of $I$, but it is not a cell complex, as an edge is missing. Using discrete Morse theory, we will be showing what the correct cell complex is to support a minimal resolution of $I$.  

\end{example}

A tool for constructing free resolutions that are smaller than the Taylor resolution is Discrete Morse Theory, as defined in~\cite{BW2002cellularresolution}. The idea is to begin with the face poset of the Taylor complex and systematically eliminate pairs of faces. In \cref{exa:polyhedral}, we apply this technique to show that the minimal free resolution of $I$ in \cref{example} is supported on a polyhedral cell complex.

%We next use Discrete Morse Theory (as defined in~\cite{BW2002cellularresolution}) to show that the minimal free resolution of $I$ in \cref{example} is supported on a polyhedral cell complex. The idea is to start from the face poset of the Taylor complex and eliminate faces in pair systematically.

%%%%%%%%%%%%%%%%%%%%%%%%%%%%%%%%%%%%%%%%%%%%%%%%%%%%%%%%%%%%%%%%%%%%%%%%%
\section{Discrete Morse theory}\label{sec:morsetheory}
%%%%%%%%%%%%%%%%%%%%%%%%%%%%%%%%%%%%%%%%%%%%%%%%%%%%%%%%%%%%%%%%%%%%%%%%%
    Starting from a cell complex $X$,
    we define a \emph{directed graph} $G_X$ on the set of cells of $X$ whose set of directed edges $E_X$ is given by $\sigma\to \sigma'$ for $\sigma'\subset \sigma$ and $\dim(\sigma) = \dim(\sigma')+1$.   
    A \emph{matching} of $G_X$ is a set $M\subseteq E_X$ with the property that each cell of $X$ occurs in at most one edge of $M$. For a matching $M$, let $G_X^M$ be the graph on cells of $X$ whose edge set is
    \[
    E_X^M = (E_X\setminus M) \cup \{\sigma'\to \sigma : \sigma\to \sigma' \in M\}.
    \]
    
    We say that $M$ is \emph{acyclic} if $G_X^M$ is acyclic i.e., it does not contain directed cycles. 
    Given an acyclic matching $M$, the cells of $X$ that do not appear in the edges of $M$ are called the \emph{$M$-critical cells}.
    For an acyclic matching $M$ of $G_X$, a directed path $\sigma_0\to \cdots \to \sigma_m$ in $G_X^M$ is called a \emph{gradient path}.

    Let $I\subseteq S$ be a monomial ideal that is minimally generated by $m_1,\ldots,m_r$. 
    Assume that $X$ is a subset of $\Taylor(I)$. We say that a matching $M$ on $G_X$ is 
    \emph{homogeneous} if for all $\sigma\to \sigma' \in M$, we have $\lcm(\sigma) = \lcm(\sigma')$, where $\lcm(\sigma) = \lcm(m_i : i\in \sigma).$ For convenience, we interpret the faces of $\Taylor(I)$ as cells of $2^{[r]}$.

\begin{theorem}[{\bf The Morse Complex $X_M$}{\cite[Proposition\ 1.2, Theorem\ 1.3]{BW2002cellularresolution}}]\label{thm:BWcellular}
    If $I$ is a monomial ideal, $X$ is the Taylor complex of $I$, and $M$ is a homogeneous acyclic matching
of $G_X$, then there is a CW complex $X_M$ which supports a multigraded free resolution of $I$.
The $i$-cells of $X_M$ are in one-to-one correspondence with the $M$-critical $i$-cells of $X$.
\end{theorem}

We will refer to the cell complex $X_M$ as the \emph{Morse complex} of the matching $M$.
\cref{thm:BWcellular} in particular indicates that the larger the matching $M$ is, the smaller the Morse complex $X_M$ will be, and in particular, the closer the cellular resolution will be to a minimal one. We are therefore interested in \emph{maximal} homogeneous acyclic matchings of the face poset of a Taylor complex.

    For an $M$-critical $i$-cell $\sigma$ of $X$, let us denote the corresponding $i$-cell of $X_M$ by $\sigma_M.$ The following lemma gives a concrete description of the ordering of cells in $X_M.$  

\begin{lemma}[Proposition~7.3, \cite{BW2002cellularresolution}]\label{lemma:subsetgradientpath}
    Let $X$ be the Taylor complex of $I$ and $M$ be a homogeneous acyclic matching of $G_X.$
    Let $\sigma$ and $\tau$ be two $M$-critical cells of $X$ such that $\dim(\tau) = \dim(\sigma)+1.$
    Then, $\sigma_M \prec \tau_M$ if and only if there is a gradient path in $G_X^M$ as
    
    $$\begin{array}{ccccccccccccc}
   \tau=\sigma_0&& && \sigma_2 &&&&&& \sigma_{m-1} &&  \\
    &\searrow&&\nearrow &&\searrow &&&&\nearrow &&\searrow & \\ 
    &&\sigma_1&&&& \sigma_3 &\ldots&\sigma_{m-2}&&&& \sigma_{m}=\sigma \\ 
    \end{array}.$$
    
\end{lemma}

\cref{lemma:subsetgradientpath} is essentially saying that the  nonsimplicial cells of the Morse complex correspond to the faces of the Taylor complex which lose a subface  through the matching $M$. This observation is summarized in the next statement which generalizes [Proposition~5.2, \cite{ChauSelvi}].

\begin{proposition}[{\bf Cells of the Morse complex}]\label{lem:2-step-path}
Let $X$ be the Taylor complex of $I$ and $M$ be a homogeneous acyclic matching of $G_X.$
Suppose $\tau\to \sigma\in M$, and   
$\omega$ is an $M$-critical face of $X$ such that $\dim(\omega)=\dim(\tau)=d$ and $\sigma \subset \omega$.
\begin{enumerate}
    \item For every $M$-critical face $\gamma \subset \tau$ with $\dim (\gamma)=d-1$, we have $\gamma_M\prec \omega_M.$
In other words, $\omega_M$ inherits, as sub-cells, all $M$-critical faces of $\tau$ of dimension $d-1$.

\item If $\tau$ and all of its faces are  $M$-critical, then for $\sigma\in X,$ we have $\sigma_M\prec \tau_M$ if and only if $\sigma \subset \tau$. 
In other words, $\tau_M$ is a simplex.
\end{enumerate}
\end{proposition}

\begin{proof}
    %It suffices to prove the case where $\dim(\sigma') = \dim(\tau)-1.$ Let $\sigma=\sigma_1\subset\cdots \subset \sigma_r=\tau'$ be a chain in $X$ such that $\dim(\sigma_{i}) = \dim(\sigma_{i+1})-1$ for $i=1,\ldots,r-1.$  
    The proof follows from \cref{lemma:subsetgradientpath} as we have the gradient path 
    %$\tau=\sigma_r \to \sigma_{r-1}\to  \cdots \to \sigma_2\to \sigma_1 =\sigma\to \tau \to \sigma'$.
$$\begin{array}{cccccccc}
    \omega &&&& \tau&\\
    &\searrow &&\nearrow &&\searrow &\\ 
    &&\sigma&&&& \gamma&\\ 
    \end{array}.$$
\end{proof}

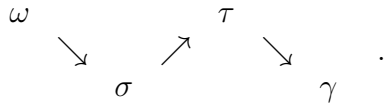
\begin{figure}[h]
\centering
\begin{tikzpicture}[scale=.82]
\node [draw, circle, fill=white, inner sep=1pt, label=above:{\tiny{{${\{123456\}}$}}}] (123456) at (9.5,5) {};

\node [draw, circle, fill=white, inner sep=1pt, label=above:{\tiny{{${\{12345\}}$}}}] (12345) at (4.5,4) {};
\node [draw, circle, fill=white, inner sep=1pt, label=above:{\tiny{{${\{12346\}}$}}}] (12346) at (6.5,4) {};
\node [draw, circle, fill=white, inner sep=1pt, label=above:{\tiny{{${\{12356\}}$}}}] (12356) at (8.5,4) {};
\node [draw, circle, fill=white, inner sep=1pt, label=above:{\tiny{{${\{12456\}}$}}}] (12456) at (10.5,4) {};
\node [draw, circle, fill=white, inner sep=1pt, label=above:{\tiny{{${\{13456\}}$}}}] (13456) at (12.5,4) {};
\node [draw, circle, fill=white, inner sep=1pt, label=above:{\tiny{{${\{23456\}}$}}}] (23456) at (14.5,4) {};

\node [draw, circle, fill=white, inner sep=1pt,label=below:{\tiny{{${\{1234\}}$}}}] (1234) at (1.5,3) {};
\node [draw, circle, fill=white, inner sep=1pt, label=below:{\tiny{{${\{1235\}}$}}}] (1235) at (2.6,3) {};
\node [draw, circle, fill=white, inner sep=1pt, label=below:{\tiny{{${\{1236\}}$}}}] (1236) at (3.7,3) {};
\node [draw, circle, fill=white, inner sep=1pt, label=below:{\tiny{{${\{1245\}}$}}}] (1245) at (4.8,3) {};
\node [draw, circle, fill=white, inner sep=1pt, label=below:{\tiny{${\{1246\}}$}}] (1246) at (5.9,3) {};
\node [draw, circle, fill=white, inner sep=1pt, label=below:{\tiny{${\{1256\}}$}}] (1256) at (7,3) {};
\node [draw, circle, fill=white, inner sep=1pt, label=below:{\tiny{${\{1345\}}$}}] (1345) at (8.1,3) {};
\node [draw, circle, fill=white, inner sep=1pt, label=below:{\tiny{${\{1346\}}$}}] (1346) at (9.2,3) {};
\node [draw, circle, fill=white, inner sep=1pt, label=below:{\tiny{{${\{1356\}}$}}}] (1356) at (10.3,3) {};
\node [draw, circle, fill=white, inner sep=1pt, label=below:{\tiny{${\{1456\}}$}}] (1456) at (11.4,3) {};
\node [draw, circle, fill=white, inner sep=1pt, label=below:{\tiny{{${\{2345}\}$}}}] (2345) at (12.5,3) {};
\node [draw, circle, fill=white, inner sep=1pt, label=below:{\tiny{${\{2346}\}$}}] (2346) at (13.6,3) {};
\node [draw, circle, fill=white, inner sep=1pt, label=below:{\tiny{{${\{2356}\}$}}}] (2356) at (14.7,3) {};
\node [draw, circle, fill=white, inner sep=1pt, label=below:{\tiny{{${\{2456}\}$}}}] (2456) at (15.8,3) {};
\node [draw, circle, fill=white, inner sep=1pt, label=below:{\tiny{{${\{3456}\}$}}}] (3456) at (16.9,3) {};

\node [draw, circle, fill=white, inner sep=1pt, label=below:{\tiny{{${\{123\}}$}}}] (123) at (0,2) {};
\node [draw, circle, fill=white, inner sep=1pt, label=below:{\tiny{${\{124\}}$}}] (124) at (1,2) {};
\node [draw, circle, fill=white, inner sep=1pt, label=below:{\tiny{${\{125\}}$}}] (125) at (2,2) {};
\node [draw, circle, fill=white, inner sep=1pt, label=below:{\tiny{${\{126\}}$}}] (126) at (3,2) {};
\node [draw, circle, fill=white, inner sep=1pt, label=below:{\tiny{${\{134\}}$}}] (134) at (4,2) {};
\node [draw, circle, fill=white, inner sep=1pt, label=below:{\tiny{${\{135\}}$}}] (135) at (5,2) {};
\node [draw, circle, fill=white, inner sep=1pt, label=below:{\tiny{${\{136\}}$}}] (136) at (6,2) {};
\node [draw, circle, fill=white, inner sep=1pt, label=below:{\tiny{${\{145\}}$}}] (145) at (7,2) {};
\node [draw, circle, fill=white, inner sep=1pt, label=below:{\tiny{${\{146\}}$}}] (146) at (8,2) {};
\node [draw, circle, fill=white, inner sep=1pt, label=below:{\tiny{${\{156\}}$}}] (156) at (9,2) {};
\node [draw, circle, fill=white, inner sep=1pt, label=below:{\tiny{${\{234\}}$}}] (234) at (10,2) {};
\node [draw, circle, fill=white, inner sep=1pt, label=below:{\tiny{${\{235\}}$}}] (235) at (11,2) {};
\node [draw, circle, fill=white, inner sep=1pt, label=below:{\tiny{${\{236\}}$}}] (236) at (12,2) {};
\node [draw, circle, fill=white, inner sep=1pt,label=below:{\tiny{{${\{245\}}$}}}] (245) at (13,2) {};
\node [draw, circle, fill=white, inner sep=1pt,label=below:{\tiny{${\{246\}}$}}] (246) at (14,2) {};
\node [draw, circle, fill=white, inner sep=1pt,label=below:{\tiny{${\{256\}}$}}] (256) at (15,2) {};
\node [draw, circle, fill=white, inner sep=1pt,label=below:{\tiny{${\{345\}}$}}] (345) at (16,2) {};
\node [draw, circle, fill=white, inner sep=1pt,label=below:{\tiny{${\{346\}}$}}] (346) at (17,2) {};
\node [draw, circle, fill=white, inner sep=1pt,label=below:{\tiny{{${\{356\}}$}}}] (356) at (18,2) {};
\node [draw, circle, fill=white, inner sep=1pt,label=below:{\tiny{${\{456\}}$}}] (456) at (19,2) {};

\node [draw, circle, fill=black, inner sep=1pt,label=below:{\tiny{${\{12\}}$}}] (12) at (1.5,1) {};
\node [draw, circle, fill=white, inner sep=1pt, label=below:{\tiny{${\{13\}}$}}] (13) at (2.6,1) {};
\node [draw, circle, fill=white, inner sep=1pt, label=below:{\tiny{${\{14\}}$}}] (14) at (3.7,1) {};
\node [draw, circle, fill=white, inner sep=1pt, label=below:{\tiny{${\{15\}}$}}] (15) at (4.8,1) {};
\node [draw, circle, fill=black, inner sep=1pt, label=below:{\tiny{${\{16\}}$}}] (16) at (5.9,1) {};
\node [draw, circle, fill=black, inner sep=1pt, label=below:{\tiny{${\{23\}}$}}] (23) at (7,1) {};
\node [draw, circle, fill=white, inner sep=1pt, label=below:{\tiny{${\{24\}}$}}] (24) at (8.1,1) {};
\node [draw, circle, fill=white, inner sep=1pt, label=below:{\tiny{${\{25\}}$}}] (25) at (9.2,1) {};
\node [draw, circle, fill=white, inner sep=1pt, label=below:{\tiny{${\{26\}}$}}] (26) at (10.3,1) {};
\node [draw, circle, fill=black, inner sep=1pt, label=below:{\tiny{${\{34\}}$}}] (34) at (11.4,1) {};
\node [draw, circle, fill=white, inner sep=1pt, label=below:{\tiny{${\{35}\}$}}] (35) at (12.5,1) {};
\node [draw, circle, fill=white, inner sep=1pt, label=below:{\tiny{${\{36}\}$}}] (36) at (13.6,1) {};
\node [draw, circle, fill=black, inner sep=1pt, label=below:{\tiny{${\{45}\}$}}] (45) at (14.7,1) {};
\node [draw, circle, fill=white, inner sep=1pt, label=below:{\tiny{${\{46}\}$}}] (46) at (15.8,1) {};
\node [draw, circle, fill=black, inner sep=1pt, label=below:{\tiny{${\{56}\}$}}] (56) at (16.9,1) {};

\node [draw, circle, fill=black, inner sep=1pt, label=below:{\tiny{${\{1\}}$}}] (1) at (4.5,0) {};
\node [draw, circle, fill=black, inner sep=1pt, label=below:{\tiny{${\{2\}}$}}] (2) at (6.5,0) {};
\node [draw, circle, fill=black, inner sep=1pt, label=below:{\tiny{${\{3\}}$}}] (3) at (8.5,0) {};
\node [draw, circle, fill=black, inner sep=1pt, label=below:{\tiny{${\{4\}}$}}] (4) at (10.5,0) {};
\node [draw, circle, fill=black, inner sep=1pt, label=below:{\tiny{${\{5\}}$}}] (5) at (12.5,0) {};
\node [draw, circle, fill=black, inner sep=1pt, label=below:{\tiny{${\{6\}}$}}] (6) at (14.5,0) {};

\draw [color=red,line width=1pt, <-] (13)--(123); 
\draw [color=red,line width=1pt, <-] (14)--(124); 
\draw [color=red,line width=1pt, <-] (15)--(156); 
\draw [color=red,line width=1pt, <-] (25)--(245); 
\draw [color=red,line width=1pt, <-] (35)--(345); 
\draw [color=red,line width=1pt, <-] (36)--(236); 
%\draw [color=red,line width=1pt, <-] (46)--(246);

\draw [color=red,line width=1pt, <-] (125)--(1256);
\draw [color=red,line width=1pt, <-] (134)--(1234);
\draw [color=red,line width=1pt, <-] (135)--(1235); 
\draw [color=red,line width=1pt, <-] (136)--(1236); 
\draw [color=red,line width=1pt, <-] (145)--(1456); 
\draw [color=red,line width=1pt, <-] (146)--(1246); 
\draw [color=red,line width=1pt, <-] (235)--(2345);
\draw [color=red,line width=1pt, <-] (256)--(2456);
\draw [color=red,line width=1pt, <-] (346)--(2346);
\draw [color=red,line width=1pt, <-] (356)--(3456);

\draw [color=red,line width=1pt, <-] (1245)--(12456);
\draw [color=red,line width=1pt, <-] (1345)--(12345);
\draw [color=red,line width=1pt, <-] (1346)--(12346);
\draw [color=red,line width=1pt, <-] (1356)--(12356); 
\draw [color=red,line width=1pt, <-] (2356)--(23456);
\draw [color=red,line width=1pt, <-] (13456)--(123456);

\end{tikzpicture}
\caption{Homogeneous acyclic matching $M$. The vertices corresponding to the faces of $\Scarf(I)$ are filled with black.}\label{fig:matchingexample}
\end{figure}
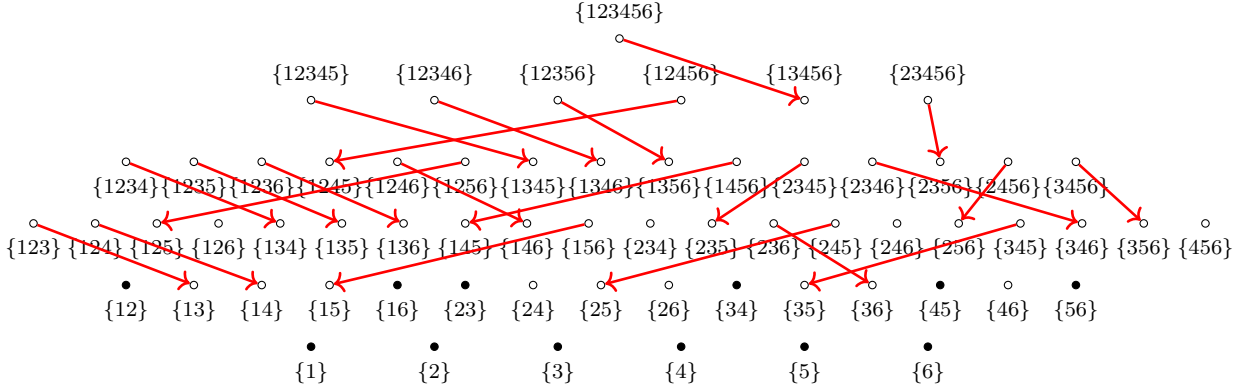

\begin{example}\label{exa:polyhedral}
    Let $I=(x^2,xy,y^2,yz,z^2,xz)$ and $\Delta$ be as in \cref{example}. Consider the homogeneous acyclic matching $M$ shown in \Cref{fig:matchingexample}. By \cref{lem:2-step-path}, $X_M$ is the simplicial complex $\Delta$ because for every $M$-critical cell of $\Taylor(I),$ all its sub-cells are also $M$-critical. 
    
    Next, let us define $M_1 = M\cup \{\{2,4,6\} \to \{4,6\}\}.$  Then the face $\{4,5,6\}$ of $\Taylor(I)$ is $M_1$-critical, but loses its subface $\{4,6\}$ through the matching $M_1$. As a result, $\{4,5,6\}_{M_1}$ is a nonsimplicial face of $X_{M_1}$, and has as subcells all $\sigma_{M_1}$ which are $M_1$-critical subfaces of $\{4,5,6\}$ or $\{2, 4,6\}$, namely: $\{2,4\}$, $\{2,6\}$, $\{4,5\}$, $\{5,6\}$ and all their vertices. The polyhedral cell complex  $X_{M_1}$ is given below on the right by \cref{lem:2-step-path}.
    
    Moreover, $M_1$ is a maximal homogeneous acyclic matching which by \cref{thm:BWcellular} supports a free resolution of $I$. Macaulay2~\cite{M2} computations show that it is the minimal free resolution of $I$.
\end{example}
%\begin{figure}[h]
%    \centering
%    \includegraphics[width=0.9\linewidth]{Screenshot 2024-10-28 000556.png}
%    \caption{Enter Caption}
%    \label{fig:enter-label}
%\end{figure}
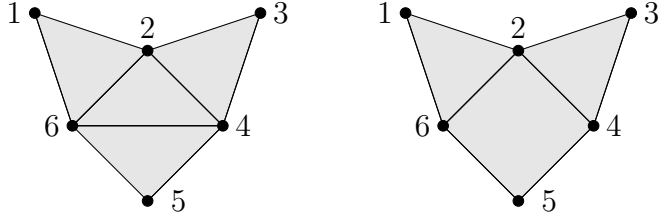
\[
\begin{tabular}{lr}
        \begin{tikzpicture}[scale = 1]
\tikzstyle{point}=[inner sep=0pt]
\node (a)[point,label=above:$2$] at (0,1) {}; 
\node (b)[point,label=left:$6$] at (-1,0) {}; 
\node (c)[point,label=right:$4$] at (1,0) {}; 
\node (d)[point,label=left:$1$] at (-1.5,1.5) {}; 
\node (e)[point,label=right:$3$] at (1.5,1.5) {}; 
\node (f)[point,label=right:$\ 5$] at (0,-1) {}; 
\draw [fill=gray!20](a.center) -- (b.center) -- (c.center);
\draw [fill=gray!20](a.center) -- (b.center) -- (d.center);
\draw [fill=gray!20](a.center) -- (c.center) -- (e.center);
\draw [fill=gray!20](b.center) -- (c.center) -- (f.center);
\draw (a.center) -- (b.center);
\draw (a.center) -- (c.center);
\draw (a.center) -- (d.center);
\draw (a.center) -- (e.center);
\draw (b.center) -- (c.center);
\draw (b.center) -- (d.center);
\draw (b.center) -- (f.center);
\draw (c.center) -- (e.center);
\draw (c.center) -- (f.center);
% \node [label=$xyz$] at (0,0) {}; 
% \node [label=$x^2yz$] at (-.8,.4) {}; 
% \node [label=$xy^2z$] at (0.8,.4) {}; 
% \node [label=$xyz^2$] at (0,-.8) {}; 
\foreach \p in {a,b,c,d,e,f}  \filldraw[black] (\p) circle (2pt);
        \end{tikzpicture}
\quad & \quad 

       \begin{tikzpicture}[scale = 1]
\tikzstyle{point}=[inner sep=0pt]
\node (a)[point,label=above:$2$] at (0,1) {}; 
\node (b)[point,label=left:$6$] at (-1,0) {}; 
\node (c)[point,label=right:$4$] at (1,0) {}; 
\node (d)[point,label=left:$1$] at (-1.5,1.5) {}; 
\node (e)[point,label=right:$3$] at (1.5,1.5) {}; 
\node (f)[point,label=right:$\ 5$] at (0,-1) {}; 
\draw [fill=gray!20](a.center) -- (b.center) -- (f.center) --(c.center);
\draw [fill=gray!20](a.center) -- (b.center) -- (d.center);
\draw [fill=gray!20](a.center) -- (c.center) -- (e.center);
%\draw [fill=gray!20](b.center) -- (c.center) -- (f.center);
\draw (a.center) -- (b.center);
\draw (a.center) -- (c.center);
\draw (a.center) -- (d.center);
\draw (a.center) -- (e.center);
%\draw (b.center) -- (c.center);
\draw (b.center) -- (d.center);
\draw (b.center) -- (f.center);
\draw (c.center) -- (e.center);
\draw (c.center) -- (f.center);
\foreach \p in {a,b,c,d,e,f}
          \filldraw[black] (\p) circle (2pt);
        \end{tikzpicture}
\end{tabular}
\]

%%%%%%%%%%%%%%%%%%%%%%%%%%%%%%%%%%%%%%%%%%%%%%%%%%%%%%%%%%%%%%%%%%%%%%%%%
\section{When are Morse complexes polyhedral?}\label{sec:4generators}
%%%%%%%%%%%%%%%%%%%%%%%%%%%%%%%%%%%%%%%%%%%%%%%%%%%%%%%%%%%%%%%%%%%%%%%%%

Given the algorithmic nature of the construction of the Morse complex, which is done via pairwise elimination of face and facets from the Taylor simplex, a natural question is: is there always a way to do the matchings so that the Morse complex is polyhedral? The answer in general is negative as we show via an example with $6$ generators in this section (\cref{nonexample}), and positive for ideals with at most $4$ monomial generators (\cref{4gens}). We note that it was already shown in~\cite{Roya} that for monomial ideals with up to $4$ generators and all their monomial Artinian reductions have minimal cellular Morse resolutions. Recently, Montaner et al.~\cite{montaner2025morse} and Chau~\cite{chau2025minimal} demonstrated that ideals with $5$ monomial generators and all their monomial Artinian reductions have cellular minimal Morse resolutions.

Let $I$ be an ideal in a polynomial ring $S=K[x_1,\ldots,x_n]$ minimally generated by $r$ monomials.
Let $X$ denote the Taylor complex of $I$.
%We denote $2^{[r]}$ the set of subsets of $[r].$
We say that the set $\{(\sigma_1,\tau_1),\ldots,(\sigma_m,\tau_m)\}$ in $X\times X$ consists of \emph{minimal homogeneous pairs} of $X$ if the following conditions are satisfied:
\begin{enumerate}
    \item 
        for each $1\leq i\leq m,$ we have $\sigma_i\subset \tau_i,$ $\dim(\tau_i) = \dim(\sigma_i)+1,$ and $\lcm(\sigma_i)=\lcm(\tau_i);$
    \item 
        for every $1\leq i\leq m,$ $\lcm(\sigma_i\setminus A)\neq \lcm(\tau_i\setminus A)$ for any non-empty subset $A\subset\sigma_i;$ 
    \item
        if there exist $\sigma,\tau\in X$ such that $\sigma\subset \tau,$ $\dim(\tau) = \dim(\sigma)+1,$ and $\lcm(\sigma)=\lcm(\tau),$ then there exists an $i$ such that $(\sigma,\tau) = (\sigma_i\cup A, \tau_i \cup A)$ for some $A \subset [r]\setminus \tau_i.$
\end{enumerate}

\begin{theorem}\label{4gens}
    Let $I$ be a monomial ideal minimally generated by $\leq 4$ generators, and $X$ be its Taylor complex. Then, there is a maximal homogeneous acyclic matching $M$ of $G_X$ such that $X_M$ is a polyhedral cell complex.
\end{theorem}

\begin{proof}
 Assume first that $I$ is minimally generated by at most $3$ generators.  Then $X$ is a triangle, an edge or a vertex. By the minimality of the generators, we know any homogeneous matching does not include the vertices of $X$, and so the only possible matching would consist of the triangle matched with an edge of $X$, in which case the resulting Morse complex is a graph, which is simplicial.  

Now assume $I$ is minimally generated by $4$ elements. 
Let $\{(\sigma_1,\tau_1),\ldots,(\sigma_m,\tau_m)\}$ be the set of minimal homogeneous pairs of $X.$ If $\tau_i =\{1,2,3,4\}$ for all $1\leq i\leq m,$ then $M = \{\{1,2,3,4\} \to \sigma_1 =\{i_1,i_2,i_3\}\}$ is a maximal homogeneous acyclic matching because all other possible matchings necessarily involve the face $\{1,2,3,4\}$. Consequently, $X_M$ is the simplicial complex on the right in \Cref{f:4-gens}, and thus we are done. 

\begin{figure}[hbt!]
\begin{tabular}{lcr}
\begin{tikzpicture}
             \filldraw[fill = black!15!white]
            (0,0) -- (2,2) -- (4,0) -- (0,0);% -- cycle;
            \draw (2,1) -- (2,2)  (2,1) -- (0,0) (2,1)--(4,0);
             
                \filldraw[] (2,2) circle (2pt) node[anchor=south]  {$i_1$};
                \filldraw[] (0,0) circle (2pt) node[anchor=north]  {$i_2$};
                \filldraw[] (4,0) circle (2pt) node[anchor=north]  {$i_3$};
                \filldraw[] (2,1) circle (2pt) node[anchor=north] {};
	\end{tikzpicture}
 \quad & \quad 
             \begin{tikzpicture}[scale=1]
              \filldraw[fill = black!15!white]
                (0,0) -- (2,2)--(2,1) -- (0,0) 
                (2,1)--(4,0) -- (0,0)-- (2,1);
                \filldraw[] (2,2) circle (2pt) node[anchor=south]  {$i_1$};
                \filldraw[] (0,0) circle (2pt) node[anchor=north]  {$i_4$};
                \filldraw[] (4,0) circle (2pt) node[anchor=north]  {$i_2$};
                \filldraw[] (2,1) circle (2pt) node[anchor=south west]  {$i_3$};
	\end{tikzpicture}
 \end{tabular}
 \caption{}\label{f:4-gens}
\end{figure}

Now, without loss of generality, assume that $\tau_1\neq \{1,2,3,4\}$. Since the $0$-cells correspond to minimal generators of $I$, $\sigma_i$ is not a 0-cell for all $1\leq i\leq m.$ So $(\sigma_1,\tau_1) = (\{i_1,i_2\},\{i_1,i_2,i_3\}).$ Consider the homogeneous acyclic matching 
$$M = \{ \tau_1 \to \sigma_1, \{1,2,3,4\} \to \sigma_1\cup \{i_4\}\} 
\qwhere 
i_4 \notin \{i_1,i_2,i_3\}.$$ 
Then $X_M$ is the simplicial complex on the right in \Cref{f:4-gens}.

If $M$ is not a maximal matching, we will consider enlargements of it. In other words, 
we show that every homogeneous acyclic matching of the simplicial complex on the right gives rise to a polyhedral cell complex. Let $M_1$ be a homogeneous acyclic matching of the simplicial complex on the right. Then $M_1$ is one of the following sets
$$\{\{i_1,i_3,i_4\} \to \sigma\}, \quad 
\{\{i_2,i_3,i_4\} \to \sigma\},    \qor 
\{\{i_1,i_3,i_4\} \to \sigma, \{i_2,i_3,i_4\} \to \sigma'\}$$
for some $\sigma$ and $\sigma'$ with $\sigma\neq \sigma'$. 
We proceed in all cases:

{\bf Case 1:} When $M_1 = \{\{i_1,i_3,i_4\} \to \sigma\}$ for some $\sigma.$ We divide this case in two following subcases. In the first subcase, if $\sigma = \{i_3,i_4\},$ then $X_{M_1}$ is the polyhedral cell complex below on the left by \Cref{lem:2-step-path}, whereas in the second subcase, if $\sigma \neq \{i_3,i_4\},$ then $X_{M_1}$ is isomorphic to the simplicial complex below on the right.

    \[	
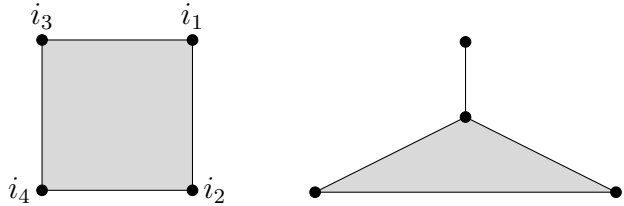
\begin{tabular}{lcr}
\begin{tikzpicture}[scale=1]
             \filldraw[fill = black!15!white]
            (0,0) -- (2,0)--(2,2) -- (0,2) -- (0,0)
            ;
            
               \filldraw[] (2,2) circle (2pt) node[anchor=south]  {$i_1$};
                \filldraw[] (0,0) circle (2pt) node[anchor=east]  {$i_4$};
                \filldraw[] (2,0) circle (2pt) node[anchor=west]  {$i_2$};
                \filldraw[] (0,2) circle (2pt) node[anchor=south]  {$i_3$};
	\end{tikzpicture}
 \quad & \quad 
             \begin{tikzpicture}[scale=1]
              \filldraw[fill = black!15!white]
                (0,0) --(2,1) --(4,0)-- (0,0) 
                ;
                \draw[] (2,2)--(2,1);
                
                 \filldraw[] (2,2) circle (2pt) node[anchor=south]  {};
                 \filldraw[] (0,0) circle (2pt) node[anchor=north]  {};
                 \filldraw[] (4,0) circle (2pt) node[anchor=north]  {};
                 \filldraw[] (2,1) circle (2pt) node[anchor=south west]  {};
	\end{tikzpicture}
 \end{tabular}
\]  

{\bf Case 2:} When $M_1 = \{\{i_2,i_3,i_4\} \to \sigma\}$ for some $\sigma.$ Use the same argument of Case 1 to get that $X_{M_1}$ is the polyhedral cell complex.

{\bf Case 3:} When $M_1 = \{\{i_1,i_3,i_4\} \to \sigma, \{i_2,i_3,i_4\} \to \sigma'\}$ for some $\sigma$ and $\sigma'$ with $\sigma\neq \sigma'$. We also divide this case in the following two subcases. In the first subcase, if either $\sigma = \{i_3,i_4\}$ or $\sigma' = \{i_3,i_4\},$ then $X_{M_1}$ is the simplicial complex below on the left. In the second subcase, if both $\sigma \neq \{i_3,i_4\}$ and $\sigma' \neq \{i_3,i_4\},$ then $X_{M_1}$ is the simplicial complex below on the right.
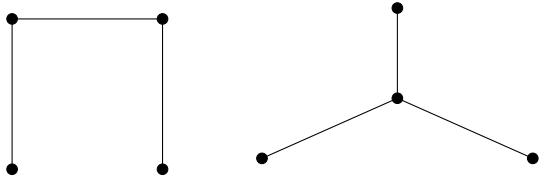
\begin{center}
\begin{tabular}{lcr}
\begin{tikzpicture}[scale=1]
             \draw[]
               (0,0) -- (0,2)--(2,2) -- (2,0)
            ;
            
                \filldraw[] (2,2) circle (2pt) node[anchor=south]  {};
                \filldraw[] (0,0) circle (2pt) node[anchor=east]  {};
                \filldraw[] (2,0) circle (2pt) node[anchor=west]  {};
                \filldraw[] (0,2) circle (2pt) node[anchor=south]  {};
	\end{tikzpicture}
 \quad & \quad 
             \begin{tikzpicture}[scale=1]
                \draw[] (2,2)--(2,.8) -- (.2,0)  (2,.8)--(3.8,0);
                
                \filldraw[] (2,2) circle (2pt) node[anchor=south]  {};
                \filldraw[] (.2,0) circle (2pt) node[anchor=north]  {};
                \filldraw[] (3.8,0) circle (2pt) node[anchor=north]  {};
                \filldraw[] (2,.8) circle (2pt) node[anchor=south west]  {};
	\end{tikzpicture}
 \end{tabular}
 \end{center}
\end{proof}

%%%%%%%%%%%%%%%%%%%%%%%%%%%%%%%%%%%%%%%%%%%%%%%%%%%%%%%%%%%%%%%%%%%%%%%%%
\section{An ideal with a no polyhedral minimal Morse resolution}\label{sec:non-polyhedralexample}
%%%%%%%%%%%%%%%%%%%%%%%%%%%%%%%%%%%%%%%%%%%%%%%%%%%%%%%%%%%%%%%%%%%%%%%%%

Consider the monomial ideal $I = (m_1,\ldots,m_6)\subset S = K[x_1,\ldots,x_{12}],$ where 
\begin{eqnarray}\label{e:gens}
m_1 = x_3x_4x_5x_6x_7, &
m_2 = x_2x_3x_{10}x_{11}, &
m_3 = x_1x_6x_9, \\ \nonumber
m_4 = x_1x_2x_4x_5x_6x_{10}, &
m_5 = x_4x_7x_8x_{10}, &
m_6 = x_2x_5x_{12}.
\end{eqnarray}

It is known that every multigraded free resolution of $I$ must contain the Scarf multidegrees of $I$.  This means, in particular, that 
every cell complex supporting a free resolution of $I$ must contain the Scarf (simplicial) complex of $I$. 
The facets of the simplicial complex $\Scarf(I)$ are as follows:
\begin{equation}\label{e:Scarf-facets}
\{1, 4, 5, 6\}, \{2, 3, 5\}, \{1, 2, 4, 6\}, \{1, 3, 4, 5\}, \{2, 3, 4, 6\}, \{1, 3, 4, 6\}, \{1, 2, 5, 6\}.
\end{equation}

\begin{figure}[hbt!]
    \centering
    \begin{subfigure}[t]{0.32\textwidth}
        \centering
        \begin{tikzpicture}[scale=1]
    \draw (0,0) -- (1,2)
                -- (2,0)
                -- (1.2,-0.5)
                -- (0,0);

    \draw (1,2) -- (1.2,-0.5);

    \draw[densely dotted] (0,0) -- (2,0);

    \draw[dashed] (0,0) -- (-1,2)
                        -- (1.2,-0.5);

    \draw[dashed] (-1,2) -- (2,0);

    \draw[dashed] (-1,2) -- (1,2);

    \filldraw[] (0,0) circle (2pt)
                node[anchor=south]{$1$};
    \filldraw[] (1.2,-0.5) circle (2pt)
                node[anchor=west]{$4$};
    \filldraw[] (2,0) circle (2pt)
                node[anchor=south]{$6$};
    \filldraw[] (1,2) circle (2pt)
                node[anchor=south]{$5$};
    \filldraw[] (-1,2) circle (2pt)
                node[anchor=south]{$2$};
\end{tikzpicture}
        \caption{Self-intersection in $\mathbb{R}^3$.}
        \label{fig:sub:intersect}
    \end{subfigure}
    \begin{subfigure}[t]{0.32\textwidth}
        \centering
        \begin{tikzpicture}[scale=1]
    \filldraw[fill = black!15!white]
             (1,1) -- (3,1)
                   -- (2,0)
                   -- (1,1);

    \draw (0,0) -- (1,1)
                -- (2,2)
                -- (3,1)
                -- (1,1)
                -- (2,0)
                -- (3,1)
                -- (4,0)
                -- (0,0);

    \filldraw[] (0,0) circle (2pt)
                node[anchor=north]{$1$};
    \filldraw[] (1,1) circle (2pt)
                node[anchor=south]{$2$};
    \filldraw[] (2,0) circle (2pt)
                node[anchor=north]{$3$};
    \filldraw[] (2,2) circle (2pt)
            node[anchor=south]{$4$};
    \filldraw[] (3,1) circle (2pt)
                node[anchor=south]{$5$};
    \filldraw[] (4,0) circle (2pt)
                node[anchor=north]{$6$};
\end{tikzpicture}
        \caption{The missing $3$-cells.}
        \label{fig:sub:base}
    \end{subfigure}
    \begin{subfigure}[t]{0.32\textwidth}
        \centering
        \begin{tikzpicture}[scale=1]
    \draw (0,0) -- (1.2,-0.5)
                -- (2,0);

    \draw[densely dotted] (0,0) -- (2,0);

    \draw (0,0) -- (1,2)
                -- (2,0)
                -- (1,-2)
                -- (0,0);

    \draw (1,2) -- (1.2,-0.5)
                -- (1,-2);

    \draw[dashed] (1,2) -- (3,-0.5)
                        -- (1,-2)
                        -- (1,2);

    \draw[dashed] (2,0) -- (3,-0.5);

    \draw[dashed, ultra thick] (1,2) -- (2,0)
                                     -- (1,-2);

    \filldraw[] (0,0) circle (2pt)
                node[anchor=south]{$3$};
    \filldraw[] (2,0) circle (2pt)
                node[anchor=south]{$1$};
    \filldraw[] (1.2,-0.5) circle (2pt)
                node[anchor=west]{$2$};
    \filldraw[] (1,2) circle (2pt)
                node[anchor=south]{$5$};
    \filldraw[] (1,-2) circle (2pt)
                node[anchor=north]{$4$};
    \filldraw[] (3,-0.5) circle (2pt)
                node[anchor=south]{$6$};
\end{tikzpicture}

        \caption{Non-polyhedral gluing.}
        \label{fig:sub:nonpoly}
    \end{subfigure}
    \caption{Illustrative sub-complexes of $\Scarf( I )$.}
\end{figure}
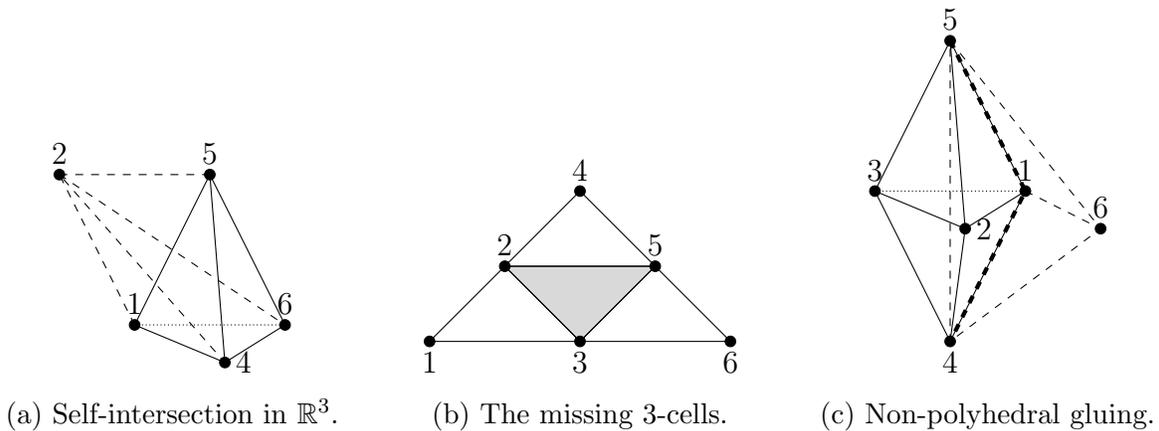

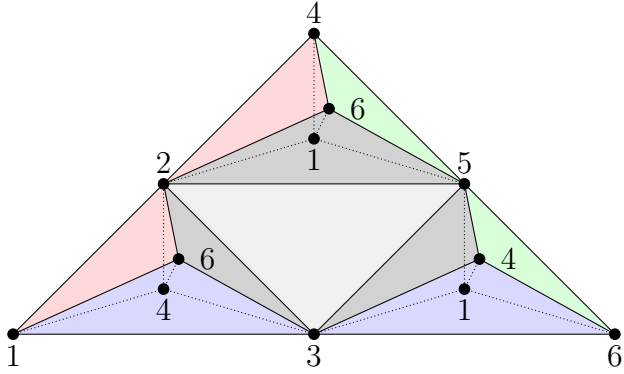
\begin{figure}[hbt!]
    \centering
    \begin{tikzpicture}[scale=1]
        % Highlighting
        \filldraw[draw = none, fill = green!15!white]
                 (8,0) -- (6.2,1)
                       -- (6,2);
        \filldraw[draw = none, fill = green!15!white]
                 (6,2) -- (4.2,3)
                       -- (4,4);

        \filldraw[draw = none, fill = red!15!white]
                 (0,0) -- (2.2,1)
                       -- (2,2);
        \filldraw[draw = none, fill = red!15!white]
                 (2,2) -- (4.2,3)
                       -- (4.2,3)
                       -- (4,4);

        \filldraw[draw = none, fill = blue!15!white]
                 (0,0) -- (2.2,1)
                       -- (4,0);
        \filldraw[draw = none, fill = blue!15!white]
                 (4,0) -- (6.2,1)
                       -- (8,0);

        \filldraw[draw = none, fill = gray!35!white]
                 (2.2,1) -- (2,2)
                         -- (4,0);
        \filldraw[draw = none, fill = gray!35!white]
                 (4,0) -- (6.2,1)
                       -- (6,2);
        \filldraw[draw = none, fill = gray!35!white]
                 (2,2) -- (6,2)
                       -- (4.2,3);

        \filldraw[draw = none, fill = gray!12!white]
                 (2,2) -- (4,0)
                       -- (6,2);

        % Framing
        \draw (0,0) -- (2,2)
                    -- (4,4)
                    -- (6,2)
                    -- (2,2)
                    -- (4,0)
                    -- (6,2)
                    -- (8,0)
                    -- (0,0);

        \filldraw[] (0,0) circle (2pt)
                    node[anchor=north]{$1$};
        \filldraw[] (2,2) circle (2pt)
                    node[anchor=south]{$2$};
        \filldraw[] (4,0) circle (2pt)
                    node[anchor=north]{$3$};
        \filldraw[] (4,4) circle (2pt)
                    node[anchor=south]{$4$};
        \filldraw[] (6,2) circle (2pt)
                    node[anchor=south]{$5$};
        \filldraw[] (8,0) circle (2pt)
                    node[anchor=north]{$6$};

        % First Section
        \draw (0,0) -- (2.2,1)
                    -- (2,2);
        \draw (4,0) -- (2.2,1);

        \draw[densely dotted] (0,0) -- (2,0.6)
                                    -- (2,2);
        \draw[densely dotted] (2.2,1) -- (2,0.6)
                                      -- (4,0);

        \filldraw[] (2.2,1) circle (2pt)
                    node[anchor=west]{$\ 6$};
        \filldraw[] (2,0.6) circle (2pt)
                    node[anchor=north]{$4$};

        % Second Section
        \draw (2,2) -- (4.2,3)
                    -- (4,4);
        \draw (6,2) -- (4.2,3);

        \draw[densely dotted] (2,2) -- (4,2.6)
                                    -- (4,4);
        \draw[densely dotted] (4.2,3) -- (4,2.6)
                                      -- (6,2);

        \filldraw[] (4.2,3) circle (2pt)
                    node[anchor=west]{$\ 6$};
        \filldraw[] (4,2.6) circle (2pt)
                    node[anchor=north]{$1$};

        % Third Section
        \draw (4,0) -- (6.2,1)
                    -- (6,2);
        \draw (8,0) -- (6.2,1);

        \draw[densely dotted] (4,0) -- (6,0.6)
                                    -- (6,2);
        \draw[densely dotted] (6.2,1) -- (6,0.6)
                                      -- (8,0);

        \filldraw[] (6.2,1) circle (2pt)
                    node[anchor=west]{$\ 4$};
        \filldraw[] (6,0.6) circle (2pt)
                    node[anchor=north]{$1$};
    \end{tikzpicture}
    \caption{A visualization of $\Scarf( I )$ in $\mathbb{R}^3$. Note that three of the tetrahedra in this diagram appear twice, to avoid self-inversion. This should be thought of as a quotient relation on the space. For ease of viewing, the three pairs of $3$-cells have been highlighted in red, green, and blue, respectively. The non-repeated $3$-cells appear in dark grey, whereas the single $2$-cell which is also a facet appears in light grey.}
    \label{fig:scarf}
\end{figure}

It turns out that $\Scarf( I )$ cannot embed into $\mathbb{R}^3$ without self-intersections.
For example, consider the sub-complex in \Cref{fig:sub:intersect}, which contains facets $\{ 1, 2, 4, 6 \}$, $\{ 1, 2, 5, 6 \}$, and $ \{ 1, 4, 5, 6 \}$.
Clearly, it is impossible to place the vertex $2$ without two of the tetrahedra intersecting.
In lieu of an embedding into $\mathbb{R}^3$, we consider the space illustrated in \Cref{fig:scarf}.
By equating the three pairs of tetrahedra in this figure, $\Scarf( I )$ is obtained.

Computations conducted using SageMath~\cite{sagemath} show that $\widetilde{H}_2(\Scarf(I),K)\cong K$; thus, $\Scarf(I)$ is not acyclic and therefore cannot  support a free resolution of $I$. In terms of its structure, the $f$-vector of $\Scarf(I)$ is $(1,6, 15, 17,6)$, whereas the Betti table of $S/I$, computed using Macaulay2~\cite{M2}, is the following:
$$
\begin{matrix}
      & 0 & 1 & 2 & 3 & 4\\
     \text{total:} & 1 & 6 & 15 & 17 & 7\\
     \hline
     0: & 1 & . & . & . & .\\
     1: & . & . & . & . & .\\
     2: & . & 2 & . & . & .\\
     3: & . & 2 & . & . & .\\
     4: & . & 1 & 2 & . & .\\
     5: & . & 1 & 9 & 1 & .\\
     6: & . & . & 4 & 15 & 6\\
     7: & . & . & . & 1 & 1
\end{matrix}
$$

In particular, the Betti vector $(1,6, 15, 17,7)$ of $S/I$ differs from the $f$-vector of $\Scarf(I)$ in only one spot. 
Thus, if the minimal free resolution of $I$ is supported on a cell complex, then one must attach exactly one $3$-cell to the Scarf complex. 
Moreover, after comparing the homogenized chain complex of $\Scarf(I)$ with the multigraded minimal free resolution of $S/I$ over $S$,  
we find that the new cell will have monomial labeling $x_1x_2\cdots x_{11}$. 

We observe that there are exactly two $3$-dimensional faces of the Taylor complex of $I$ with monomial label $x_1x_2\cdots x_{11}$; these faces are 
\begin{equation}\label{e:the-taus}
\tau=\{1,2,3,5\} 
\qand 
\tau'=\{2,3,4,5\}.
\end{equation}
Summing up our analysis above, we aim to show that there is no Morse matching $M$ on $X=\Taylor(I)$ satisfying both conditions 
(1) the set of $M$-critical faces of $X$ is $A=\Scarf(I) \cup \{\tau\}$ or $A'=\Scarf(I) \cup \{\tau'\}$, and 
(2) $X_M$ is a polyhedral cell complex.

It is worth highlighting how telling this example is: we will be showing that there is no way we can attach the $3$-cell $\tau_M$ or $\tau'_M$ to the simplicial complex $\Scarf(I)$ 
in a way that it intersects $\Scarf(I)$ in a face.

\begin{theorem}[{\bf An ideal with no polyhedral minimal Morse resolution}]\label{nonexample}
    Let $I\subset S = K[x_1,\ldots,x_{12}]$ be the monomial ideal  with generators listed in \eqref{e:gens} and let $X$ denote the Taylor complex of $I$. 
    Then, for every maximal homogeneous acyclic matching $M$ of $G_X$, the cell complex $X_M$ supports the minimal free resolution of $I$ and it is not a polyhedral cell complex.
\end{theorem}

\begin{proof}
Let 
$$\MA = \{
(\sigma,\tau) \colon  
\sigma \subset \tau \in X, \ 
\dim(\tau)=\dim(\sigma) + 1,\ 
m_\sigma=m_\tau 
        \}.$$
 Let $\gamma=\{1,2,3,4,5,6\}$. The minimal homogeneous pairs of $X$, or equivalently, the minimal elements of $\MA$, are 
\begin{equation}\label{e:minimals}
\begin{array}{l}
(\sigma_1,\tau_1) = (\{1,2,3\},\{1,2,3,4\}) =(\gamma \setminus \{4,5,6\}, \gamma \setminus \{5,6\}) \\
(\sigma_2,\tau_2) = (\{2,4,5\},\{1,2,4,5\})
=(\gamma \setminus \{1,3,6\}, \gamma \setminus \{3,6\})\\
(\sigma_3,\tau_3) = (\{3,5,6\},\{3,4,5,6\})
=(\gamma \setminus \{1,2,4\}, \gamma \setminus \{1,2\})\\
(\sigma_4,\tau_4) = (\{2,3,5,6\},\{1,2,3,5,6\})
=(\gamma \setminus \{1,4\}, \gamma \setminus \{4\}).
\end{array}
\end{equation}
This in particular implies that every pair $(\sigma, \tau) \in \MA$ can be expressed as 
\begin{equation}\label{e:MA}
(\sigma,\tau)=(\sigma_i \cup A, \tau_i \cup A)
\qforsome i\in [4] \qand 
A \subseteq [r] \setminus \tau_i.
\end{equation}
 Let $M$ be  a maximal homogeneous acyclic matchings of $G_X$. Then $M$ will contain the set of directed edges
 $$M_0=\{\tau \to \sigma \colon
 (\sigma, \tau) \in \MA, \  \sigma, \tau 
 \mbox{ do not appear in any other pair in }
 \MA \}.
 $$

 \noindent{\bf Claim:} $\gamma \to \eta \in M$ for some face $\eta$ of $X$.

We have to show that $M$ contains the pair $(\eta,\gamma)$ for some $\eta$. The pairs in the set $\MA$ that include $\gamma$  are the following: 
\[
(\gamma \setminus \{4\},\gamma) 
\qand   
(\gamma \setminus \{1\}, \gamma).
\]
Using \eqref{e:MA} and \eqref{e:minimals} we see that the other pairs in $\MA$ that include the faces $\gamma \setminus \{1\}$ and $\gamma \setminus \{4\}$ are as follows:
\begin{equation}\label{e:14}
(\gamma \setminus \{1,4\}, \gamma \setminus \{4\})
\qand 
(\gamma \setminus \{1,4\},\gamma \setminus \{1\}).
\end{equation}
Since the elements of $M$ are pairwise disjoint edges, $M$ can contain at most one of the edges 
$$
\gamma \setminus \{4\}\to \gamma \setminus \{1,4\} \qand 
\gamma \setminus \{1\}\to \gamma \setminus \{1,4\},$$
and so one of $\gamma \setminus \{4\}$ or $\gamma \setminus \{1\}$ will remain to be matched with $\gamma$. 
It follows  that 
\begin{equation}\label{e:gamma}
\gamma \to \gamma \setminus \{4\} \in M
\qor  
\gamma \to \gamma \setminus \{1\}\in M.
\end{equation}
This settles our claim. 
Now observe, again using \eqref{e:MA} and \eqref{e:minimals},  that the only pairs in $\MA$ that do not involve the  face $\gamma$ and that include a face appearing in multiple pairs are those in \eqref{e:14} as well as 
\begin{equation}\label{e:146}
(\gamma \setminus \{4,6\}, \gamma \setminus \{6\})
\qand 
(\gamma \setminus \{1,6\},\gamma \setminus \{6\}).
\end{equation}

Summarizing the discussion above, we have $M$ could be one of the following sets: 
$$
\begin{array}{l}
M_1=M_0 \cup \{
\gamma  \to \gamma\setminus \{ 4\},\quad 
\gamma \setminus \{ 1\} \to \gamma\setminus \{1,4 \},\quad 
\gamma \setminus \{ 6 \} \to \gamma\setminus \{4,6 \}
\}\\
M_2=M_0 \cup\{
\gamma  \to \gamma\setminus \{ 4\},\quad 
\gamma \setminus \{ 1\} \to \gamma\setminus \{1,4 \},\quad 
\gamma \setminus \{ 6 \} \to \gamma\setminus \{1,6 \}
\}\\
M_3=M_0 \cup\{
\gamma  \to \gamma\setminus \{ 1\},\quad 
\gamma \setminus \{ 4\} \to \gamma\setminus \{1,4 \},\quad 
\gamma \setminus \{ 6 \} \to \gamma\setminus \{4,6 \}
\}\\
M_4=M_0 \cup\{
\gamma  \to \gamma\setminus \{ 1\},\quad 
\gamma \setminus \{ 4\} \to \gamma\setminus \{1,4 \},\quad 
\gamma \setminus \{ 6 \} \to \gamma\setminus \{1,6 \}
\}
\end{array}
$$
where $M_0$ (\Cref{fig:matchingM0}) consists of the following edges
$$\begin{array}{ll}
\gamma \setminus \{ 5,6\} \to \gamma\setminus \{4,5,6 \},
&
\gamma \setminus \{ 5 \} \to \gamma\setminus \{4,5 \},
\\
\gamma \setminus \{ 3,6\} \to \gamma\setminus \{1,3,6 \},
&
\gamma \setminus \{ 3 \} \to \gamma\setminus \{1,3 \},
\\
\gamma \setminus \{1, 2\} \to \gamma\setminus \{1,2,4 \},
&
\gamma \setminus \{ 2\} \to \gamma\setminus \{2,4 \}.
\end{array}
$$

\begin{figure}[h]
\centering
\begin{tikzpicture}[scale=.82]
\node [draw, circle, fill=white, inner sep=1pt, label=above:{\tiny{{${\{123456\}}$}}}] (123456) at (9.5,5) {};
\node [draw, circle, fill=white, inner sep=1pt, label=above:{\tiny{{${\{12345\}}$}}}] (12345) at (4.5,4) {};
\node [draw, circle, fill=white, inner sep=1pt, label=above:{\tiny{{${\{12346\}}$}}}] (12346) at (6.5,4) {};
\node [draw, circle, fill=white, inner sep=1pt, label=above:{\tiny{{${\{12356\}}$}}}] (12356) at (8.5,4) {};
\node [draw, circle, fill=white, inner sep=1pt, label=above:{\tiny{{${\{12456\}}$}}}] (12456) at (10.5,4) {};
\node [draw, circle, fill=white, inner sep=1pt, label=above:{\tiny{{${\{13456\}}$}}}] (13456) at (12.5,4) {};
\node [draw, circle, fill=white, inner sep=1pt, label=above:{\tiny{{${\{23456\}}$}}}] (23456) at (14.5,4) {};
\node [draw, circle, fill=white, inner sep=1pt,label=below:{\tiny{{${\{1234\}}$}}}] (1234) at (1.5,3) {};
\node [draw, circle, fill=white, inner sep=1pt, label=below:{\tiny{{${\{1235\}}$}}}] (1235) at (2.6,3) {};
\node [draw, circle, fill=white, inner sep=1pt, label=below:{\tiny{{${\{1236\}}$}}}] (1236) at (3.7,3) {};
\node [draw, circle, fill=white, inner sep=1pt, label=below:{\tiny{{${\{1245\}}$}}}] (1245) at (4.8,3) {};
\node [draw, circle, fill=black, inner sep=1pt, label=below:{\tiny{${\{1246\}}$}}] (1246) at (5.9,3) {};
\node [draw, circle, fill=black, inner sep=1pt, label=below:{\tiny{${\{1256\}}$}}] (1256) at (7,3) {};
\node [draw, circle, fill=black, inner sep=1pt, label=below:{\tiny{${\{1345\}}$}}] (1345) at (8.1,3) {};
\node [draw, circle, fill=black, inner sep=1pt, label=below:{\tiny{${\{1346\}}$}}] (1346) at (9.2,3) {};
\node [draw, circle, fill=white, inner sep=1pt, label=below:{\tiny{{${\{1356\}}$}}}] (1356) at (10.3,3) {};
\node [draw, circle, fill=black, inner sep=1pt, label=below:{\tiny{${\{1456\}}$}}] (1456) at (11.4,3) {};
\node [draw, circle, fill=white, inner sep=1pt, label=below:{\tiny{{${\{2345}\}$}}}] (2345) at (12.5,3) {};
\node [draw, circle, fill=black, inner sep=1pt, label=below:{\tiny{${\{2346}\}$}}] (2346) at (13.6,3) {};
\node [draw, circle, fill=white, inner sep=1pt, label=below:{\tiny{{${\{2356}\}$}}}] (2356) at (14.7,3) {};
\node [draw, circle, fill=white, inner sep=1pt, label=below:{\tiny{{${\{2456}\}$}}}] (2456) at (15.8,3) {};
\node [draw, circle, fill=white, inner sep=1pt, label=below:{\tiny{{${\{3456}\}$}}}] (3456) at (16.9,3) {};

\node [draw, circle, fill=white, inner sep=1pt, label=below:{\tiny{{${\{123\}}$}}}] (123) at (0,2) {};
\node [draw, circle, fill=black, inner sep=1pt, label=below:{\tiny{${\{124\}}$}}] (124) at (1,2) {};
\node [draw, circle, fill=black, inner sep=1pt, label=below:{\tiny{${\{125\}}$}}] (125) at (2,2) {};
\node [draw, circle, fill=black, inner sep=1pt, label=below:{\tiny{${\{126\}}$}}] (126) at (3,2) {};
\node [draw, circle, fill=black, inner sep=1pt, label=below:{\tiny{${\{134\}}$}}] (134) at (4,2) {};
\node [draw, circle, fill=black, inner sep=1pt, label=below:{\tiny{${\{135\}}$}}] (135) at (5,2) {};
\node [draw, circle, fill=black, inner sep=1pt, label=below:{\tiny{${\{136\}}$}}] (136) at (6,2) {};
\node [draw, circle, fill=black, inner sep=1pt, label=below:{\tiny{${\{145\}}$}}] (145) at (7,2) {};
\node [draw, circle, fill=black, inner sep=1pt, label=below:{\tiny{${\{146\}}$}}] (146) at (8,2) {};
\node [draw, circle, fill=black, inner sep=1pt, label=below:{\tiny{${\{156\}}$}}] (156) at (9,2) {};
\node [draw, circle, fill=black, inner sep=1pt, label=below:{\tiny{${\{234\}}$}}] (234) at (10,2) {};
\node [draw, circle, fill=black, inner sep=1pt, label=below:{\tiny{${\{235\}}$}}] (235) at (11,2) {};
\node [draw, circle, fill=black, inner sep=1pt, label=below:{\tiny{${\{236\}}$}}] (236) at (12,2) {};
\node [draw, circle, fill=white, inner sep=1pt,label=below:{\tiny{{${\{245\}}$}}}] (245) at (13,2) {};
\node [draw, circle, fill=black, inner sep=1pt,label=below:{\tiny{${\{246\}}$}}] (246) at (14,2) {};
\node [draw, circle, fill=black, inner sep=1pt,label=below:{\tiny{${\{256\}}$}}] (256) at (15,2) {};
\node [draw, circle, fill=black, inner sep=1pt,label=below:{\tiny{${\{345\}}$}}] (345) at (16,2) {};
\node [draw, circle, fill=black, inner sep=1pt,label=below:{\tiny{${\{346\}}$}}] (346) at (17,2) {};
\node [draw, circle, fill=white, inner sep=1pt,label=below:{\tiny{{${\{356\}}$}}}] (356) at (18,2) {};
\node [draw, circle, fill=black, inner sep=1pt,label=below:{\tiny{${\{456\}}$}}] (456) at (19,2) {};

\node [draw, circle, fill=black, inner sep=1pt,label=below:{\tiny{${\{12\}}$}}] (12) at (1.5,1) {};
\node [draw, circle, fill=black, inner sep=1pt, label=below:{\tiny{${\{13\}}$}}] (13) at (2.6,1) {};
\node [draw, circle, fill=black, inner sep=1pt, label=below:{\tiny{${\{14\}}$}}] (14) at (3.7,1) {};
\node [draw, circle, fill=black, inner sep=1pt, label=below:{\tiny{${\{15\}}$}}] (15) at (4.8,1) {};
\node [draw, circle, fill=black, inner sep=1pt, label=below:{\tiny{${\{16\}}$}}] (16) at (5.9,1) {};
\node [draw, circle, fill=black, inner sep=1pt, label=below:{\tiny{${\{23\}}$}}] (23) at (7,1) {};
\node [draw, circle, fill=black, inner sep=1pt, label=below:{\tiny{${\{24\}}$}}] (24) at (8.1,1) {};
\node [draw, circle, fill=black, inner sep=1pt, label=below:{\tiny{${\{25\}}$}}] (25) at (9.2,1) {};
\node [draw, circle, fill=black, inner sep=1pt, label=below:{\tiny{${\{26\}}$}}] (26) at (10.3,1) {};
\node [draw, circle, fill=black, inner sep=1pt, label=below:{\tiny{${\{34\}}$}}] (34) at (11.4,1) {};
\node [draw, circle, fill=black, inner sep=1pt, label=below:{\tiny{${\{35}\}$}}] (35) at (12.5,1) {};
\node [draw, circle, fill=black, inner sep=1pt, label=below:{\tiny{${\{36}\}$}}] (36) at (13.6,1) {};
\node [draw, circle, fill=black, inner sep=1pt, label=below:{\tiny{${\{45}\}$}}] (45) at (14.7,1) {};
\node [draw, circle, fill=black, inner sep=1pt, label=below:{\tiny{${\{46}\}$}}] (46) at (15.8,1) {};
\node [draw, circle, fill=black, inner sep=1pt, label=below:{\tiny{${\{56}\}$}}] (56) at (16.9,1) {};

\node [draw, circle, fill=black, inner sep=1pt, label=below:{\tiny{${\{1\}}$}}] (1) at (4.5,0) {};
\node [draw, circle, fill=black, inner sep=1pt, label=below:{\tiny{${\{2\}}$}}] (2) at (6.5,0) {};
\node [draw, circle, fill=black, inner sep=1pt, label=below:{\tiny{${\{3\}}$}}] (3) at (8.5,0) {};
\node [draw, circle, fill=black, inner sep=1pt, label=below:{\tiny{${\{4\}}$}}] (4) at (10.5,0) {};
\node [draw, circle, fill=black, inner sep=1pt, label=below:{\tiny{${\{5\}}$}}] (5) at (12.5,0) {};
\node [draw, circle, fill=black, inner sep=1pt, label=below:{\tiny{${\{6\}}$}}] (6) at (14.5,0) {};

\draw [color=blue,line width=1pt, <-] (12346)--(1236);
\draw [color=blue,line width=1pt, <-] (12456)--(2456);
\draw [color=blue,line width=1pt, <-] (13456)--(1356);
\draw [color=blue,line width=1pt, <-] (1234)--(123);
\draw [color=blue,line width=1pt, <-] (1245)--(245);
\draw [color=blue,line width=1pt, <-] (3456)--(356);

\end{tikzpicture}
\caption{The matching $M_0$, Scarf faces are highlighted with black filled dots}\label{fig:matchingM0}
\end{figure}
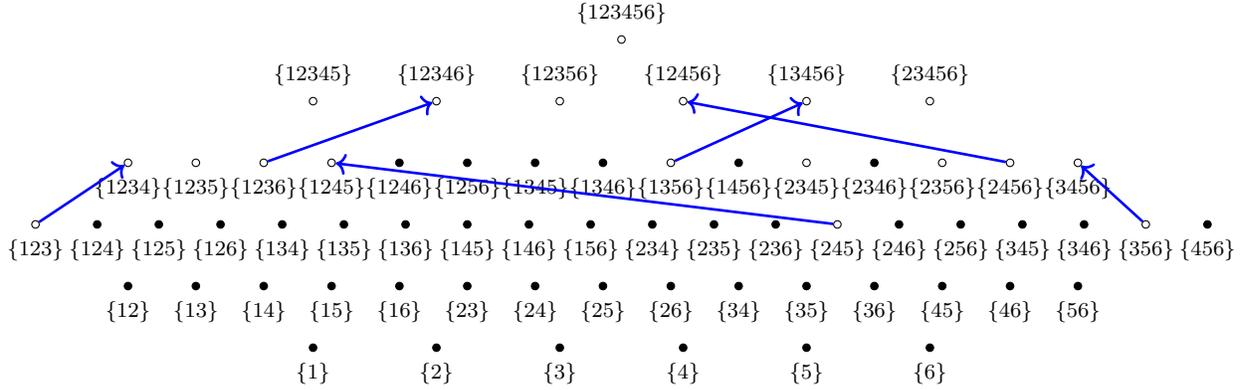

\begin{figure}[h]
\centering
\begin{tabular}{|c|c|}
\hline
\begin{tikzpicture}[scale=.82]
\node [draw, circle, fill=white, inner sep=1pt, label=above:{\tiny{{${\{123456\}}$}}}] (123456) at (5,5) {};
\node [draw, circle, fill=white, inner sep=1pt, label=above:{\tiny{{${\{12345\}}$}}}] (12345) at (1.5,4) {};
\node [draw, circle, fill=white, inner sep=1pt, label=above:{\tiny{{${\{12346\}}$}}}] (12346) at (3,4) {};
\node [draw, circle, fill=white, inner sep=1pt, label=above:{\tiny{{${\{12356\}}$}}}] (12356) at (4.5,4) {};
\node [draw, circle, fill=white, inner sep=1pt, label=above:{\tiny{{${\{12456\}}$}}}] (12456) at (6,4) {};
\node [draw, circle, fill=white, inner sep=1pt, label=above:{\tiny{{${\{13456\}}$}}}] (13456) at (7.5,4) {};
\node [draw, circle, fill=white, inner sep=1pt, label=above:{\tiny{{${\{23456\}}$}}}] (23456) at (9,4) {};
\node [draw, circle, fill=white, inner sep=1pt,label=below:{\tiny{{${\{1234\}}$}}}] (1234) at (1,3) {};
\node [draw, circle, fill=white, inner sep=1pt, label=below:{\tiny{{${\{1235\}}$}}}] (1235) at (2,3) {};
\node [draw, circle, fill=white, inner sep=1pt, label=below:{\tiny{{${\{1236\}}$}}}] (1236) at (3,3) {};
\node [draw, circle, fill=white, inner sep=1pt, label=below:{\tiny{{${\{1245\}}$}}}] (1245) at (4,3) {};
\node [draw, circle, fill=white, inner sep=1pt, label=below:{\tiny{{${\{1356\}}$}}}] (1356) at (5,3) {};
\node [draw, circle, fill=white, inner sep=1pt, label=below:{\tiny{{${\{2345}\}$}}}] (2345) at (6,3) {};
\node [draw, circle, fill=white, inner sep=1pt, label=below:{\tiny{{${\{2356}\}$}}}] (2356) at (7,3) {};
\node [draw, circle, fill=white, inner sep=1pt, label=below:{\tiny{{${\{2456}\}$}}}] (2456) at (8,3) {};
\node [draw, circle, fill=white, inner sep=1pt, label=below:{\tiny{{${\{3456}\}$}}}] (3456) at (9,3) {};
\node [draw, circle, fill=white, inner sep=1pt, label=below:{\tiny{{${\{123\}}$}}}] (123) at (3,2) {};
\node [draw, circle, fill=white, inner sep=1pt,label=below:{\tiny{{${\{245\}}$}}}] (245) at (5,2) {};
\node [draw, circle, fill=white, inner sep=1pt,label=below:{\tiny{{${\{356\}}$}}}] (356) at (7,2) {};
\draw [color=red,line width=1pt, <-] (123456)--(12356);
\draw [color=red,line width=1pt, <-] (12345)--(1235);
\draw [color=blue,line width=1pt, <-] (12346)--(1236);
\draw [color=blue,line width=1pt, <-] (12456)--(2456);
\draw [color=blue,line width=1pt, <-] (13456)--(1356);
\draw [color=red,line width=1pt, <-] (23456)--(2356);
\draw [color=blue,line width=1pt, <-] (1234)--(123);
\draw [color=blue,line width=1pt, <-] (1245)--(245);
\draw [color=blue,line width=1pt, <-] (3456)--(356);
\end{tikzpicture}
&
\begin{tikzpicture}[scale=.82]
\node [draw, circle, fill=white, inner sep=1pt, label=above:{\tiny{{${\{123456\}}$}}}] (123456) at (5,5) {};
\node [draw, circle, fill=white, inner sep=1pt, label=above:{\tiny{{${\{12345\}}$}}}] (12345) at (1.5,4) {};
\node [draw, circle, fill=white, inner sep=1pt, label=above:{\tiny{{${\{12346\}}$}}}] (12346) at (3,4) {};
\node [draw, circle, fill=white, inner sep=1pt, label=above:{\tiny{{${\{12356\}}$}}}] (12356) at (4.5,4) {};
\node [draw, circle, fill=white, inner sep=1pt, label=above:{\tiny{{${\{12456\}}$}}}] (12456) at (6,4) {};
\node [draw, circle, fill=white, inner sep=1pt, label=above:{\tiny{{${\{13456\}}$}}}] (13456) at (7.5,4) {};
\node [draw, circle, fill=white, inner sep=1pt, label=above:{\tiny{{${\{23456\}}$}}}] (23456) at (9,4) {};
\node [draw, circle, fill=white, inner sep=1pt,label=below:{\tiny{{${\{1234\}}$}}}] (1234) at (1,3) {};
\node [draw, circle, fill=white, inner sep=1pt, label=below:{\tiny{{${\{1235\}}$}}}] (1235) at (2,3) {};
\node [draw, circle, fill=white, inner sep=1pt, label=below:{\tiny{{${\{1236\}}$}}}] (1236) at (3,3) {};
\node [draw, circle, fill=white, inner sep=1pt, label=below:{\tiny{{${\{1245\}}$}}}] (1245) at (4,3) {};
\node [draw, circle, fill=white, inner sep=1pt, label=below:{\tiny{{${\{1356\}}$}}}] (1356) at (5,3) {};
\node [draw, circle, fill=white, inner sep=1pt, label=below:{\tiny{{${\{2345}\}$}}}] (2345) at (6,3) {};
\node [draw, circle, fill=white, inner sep=1pt, label=below:{\tiny{{${\{2356}\}$}}}] (2356) at (7,3) {};
\node [draw, circle, fill=white, inner sep=1pt, label=below:{\tiny{{${\{2456}\}$}}}] (2456) at (8,3) {};
\node [draw, circle, fill=white, inner sep=1pt, label=below:{\tiny{{${\{3456}\}$}}}] (3456) at (9,3) {};
\node [draw, circle, fill=white, inner sep=1pt, label=below:{\tiny{{${\{123\}}$}}}] (123) at (3,2) {};
\node [draw, circle, fill=white, inner sep=1pt,label=below:{\tiny{{${\{245\}}$}}}] (245) at (5,2) {};
\node [draw, circle, fill=white, inner sep=1pt,label=below:{\tiny{{${\{356\}}$}}}] (356) at (7,2) {};

\draw [color=red,line width=1pt, <-] (123456)--(12356);
\draw [color=red,line width=1pt, <-] (12345)--(2345);
\draw [color=blue,line width=1pt, <-] (12346)--(1236);
\draw [color=blue,line width=1pt, <-] (12456)--(2456);
\draw [color=blue,line width=1pt, <-] (13456)--(1356);
\draw [color=red,line width=1pt, <-] (23456)--(2356);
\draw [color=blue,line width=1pt, <-] (1234)--(123);
\draw [color=blue,line width=1pt, <-] (1245)--(245);
\draw [color=blue,line width=1pt, <-] (3456)--(356);
\end{tikzpicture}
\\
$M_1$& $M_2$\\
\hline
\begin{tikzpicture}[scale=.82]
\node [draw, circle, fill=white, inner sep=1pt, label=above:{\tiny{{${\{123456\}}$}}}] (123456) at (5,5) {};
\node [draw, circle, fill=white, inner sep=1pt, label=above:{\tiny{{${\{12345\}}$}}}] (12345) at (1.5,4) {};
\node [draw, circle, fill=white, inner sep=1pt, label=above:{\tiny{{${\{12346\}}$}}}] (12346) at (3,4) {};
\node [draw, circle, fill=white, inner sep=1pt, label=above:{\tiny{{${\{12356\}}$}}}] (12356) at (4.5,4) {};
\node [draw, circle, fill=white, inner sep=1pt, label=above:{\tiny{{${\{12456\}}$}}}] (12456) at (6,4) {};
\node [draw, circle, fill=white, inner sep=1pt, label=above:{\tiny{{${\{13456\}}$}}}] (13456) at (7.5,4) {};
\node [draw, circle, fill=white, inner sep=1pt, label=above:{\tiny{{${\{23456\}}$}}}] (23456) at (9,4) {};
\node [draw, circle, fill=white, inner sep=1pt,label=below:{\tiny{{${\{1234\}}$}}}] (1234) at (1,3) {};
\node [draw, circle, fill=white, inner sep=1pt, label=below:{\tiny{{${\{1235\}}$}}}] (1235) at (2,3) {};
\node [draw, circle, fill=white, inner sep=1pt, label=below:{\tiny{{${\{1236\}}$}}}] (1236) at (3,3) {};
\node [draw, circle, fill=white, inner sep=1pt, label=below:{\tiny{{${\{1245\}}$}}}] (1245) at (4,3) {};
\node [draw, circle, fill=white, inner sep=1pt, label=below:{\tiny{{${\{1356\}}$}}}] (1356) at (5,3) {};
\node [draw, circle, fill=white, inner sep=1pt, label=below:{\tiny{{${\{2345}\}$}}}] (2345) at (6,3) {};
\node [draw, circle, fill=white, inner sep=1pt, label=below:{\tiny{{${\{2356}\}$}}}] (2356) at (7,3) {};
\node [draw, circle, fill=white, inner sep=1pt, label=below:{\tiny{{${\{2456}\}$}}}] (2456) at (8,3) {};
\node [draw, circle, fill=white, inner sep=1pt, label=below:{\tiny{{${\{3456}\}$}}}] (3456) at (9,3) {};
\node [draw, circle, fill=white, inner sep=1pt, label=below:{\tiny{{${\{123\}}$}}}] (123) at (3,2) {};
\node [draw, circle, fill=white, inner sep=1pt,label=below:{\tiny{{${\{245\}}$}}}] (245) at (5,2) {};
\node [draw, circle, fill=white, inner sep=1pt,label=below:{\tiny{{${\{356\}}$}}}] (356) at (7,2) {};

\draw [color=red,line width=1pt, <-] (123456)--(23456);
\draw [color=red,line width=1pt, <-] (12345)--(1235);
\draw [color=blue,line width=1pt, <-] (12346)--(1236);
\draw [color=red,line width=1pt, <-] (12356)--(2356);
\draw [color=blue,line width=1pt, <-] (12456)--(2456);
\draw [color=blue,line width=1pt, <-] (13456)--(1356);
%\draw [color=red,line width=1pt, <-] (23456)--(2356);
\draw [color=blue,line width=1pt, <-] (1234)--(123);
\draw [color=blue,line width=1pt, <-] (1245)--(245);
\draw [color=blue,line width=1pt, <-] (3456)--(356);
\end{tikzpicture}
&
\begin{tikzpicture}[scale=.82]
\node [draw, circle, fill=white, inner sep=1pt, label=above:{\tiny{{${\{123456\}}$}}}] (123456) at (5,5) {};
\node [draw, circle, fill=white, inner sep=1pt, label=above:{\tiny{{${\{12345\}}$}}}] (12345) at (1.5,4) {};
\node [draw, circle, fill=white, inner sep=1pt, label=above:{\tiny{{${\{12346\}}$}}}] (12346) at (3,4) {};
\node [draw, circle, fill=white, inner sep=1pt, label=above:{\tiny{{${\{12356\}}$}}}] (12356) at (4.5,4) {};
\node [draw, circle, fill=white, inner sep=1pt, label=above:{\tiny{{${\{12456\}}$}}}] (12456) at (6,4) {};
\node [draw, circle, fill=white, inner sep=1pt, label=above:{\tiny{{${\{13456\}}$}}}] (13456) at (7.5,4) {};
\node [draw, circle, fill=white, inner sep=1pt, label=above:{\tiny{{${\{23456\}}$}}}] (23456) at (9,4) {};
\node [draw, circle, fill=white, inner sep=1pt,label=below:{\tiny{{${\{1234\}}$}}}] (1234) at (1,3) {};
\node [draw, circle, fill=white, inner sep=1pt, label=below:{\tiny{{${\{1235\}}$}}}] (1235) at (2,3) {};
\node [draw, circle, fill=white, inner sep=1pt, label=below:{\tiny{{${\{1236\}}$}}}] (1236) at (3,3) {};
\node [draw, circle, fill=white, inner sep=1pt, label=below:{\tiny{{${\{1245\}}$}}}] (1245) at (4,3) {};
\node [draw, circle, fill=white, inner sep=1pt, label=below:{\tiny{{${\{1356\}}$}}}] (1356) at (5,3) {};
\node [draw, circle, fill=white, inner sep=1pt, label=below:{\tiny{{${\{2345}\}$}}}] (2345) at (6,3) {};
\node [draw, circle, fill=white, inner sep=1pt, label=below:{\tiny{{${\{2356}\}$}}}] (2356) at (7,3) {};
\node [draw, circle, fill=white, inner sep=1pt, label=below:{\tiny{{${\{2456}\}$}}}] (2456) at (8,3) {};
\node [draw, circle, fill=white, inner sep=1pt, label=below:{\tiny{{${\{3456}\}$}}}] (3456) at (9,3) {};
\node [draw, circle, fill=white, inner sep=1pt, label=below:{\tiny{{${\{123\}}$}}}] (123) at (3,2) {};
\node [draw, circle, fill=white, inner sep=1pt,label=below:{\tiny{{${\{245\}}$}}}] (245) at (5,2) {};
\node [draw, circle, fill=white, inner sep=1pt,label=below:{\tiny{{${\{356\}}$}}}] (356) at (7,2) {};
\draw [color=red,line width=1pt, <-] (123456)--(23456);
\draw [color=red,line width=1pt, <-] (12345)--(2345);
\draw [color=blue,line width=1pt, <-] (12346)--(1236);
\draw [color=red,line width=1pt, <-] (12356)--(2356);
\draw [color=blue,line width=1pt, <-] (12456)--(2456);
\draw [color=blue,line width=1pt, <-] (13456)--(1356);
%\draw [color=red,line width=1pt, <-] (23456)--(2356);
\draw [color=blue,line width=1pt, <-] (1234)--(123);
\draw [color=blue,line width=1pt, <-] (1245)--(245);
\draw [color=blue,line width=1pt, <-] (3456)--(356);
\end{tikzpicture}
\\
$M_3$& $M_4$\\
\hline
\end{tabular}
\caption{The matchings $M_1,\ldots,M_4$, excluding Scarf faces, edges in  $M_0$ are blue}\label{fig:matchings}
\end{figure}
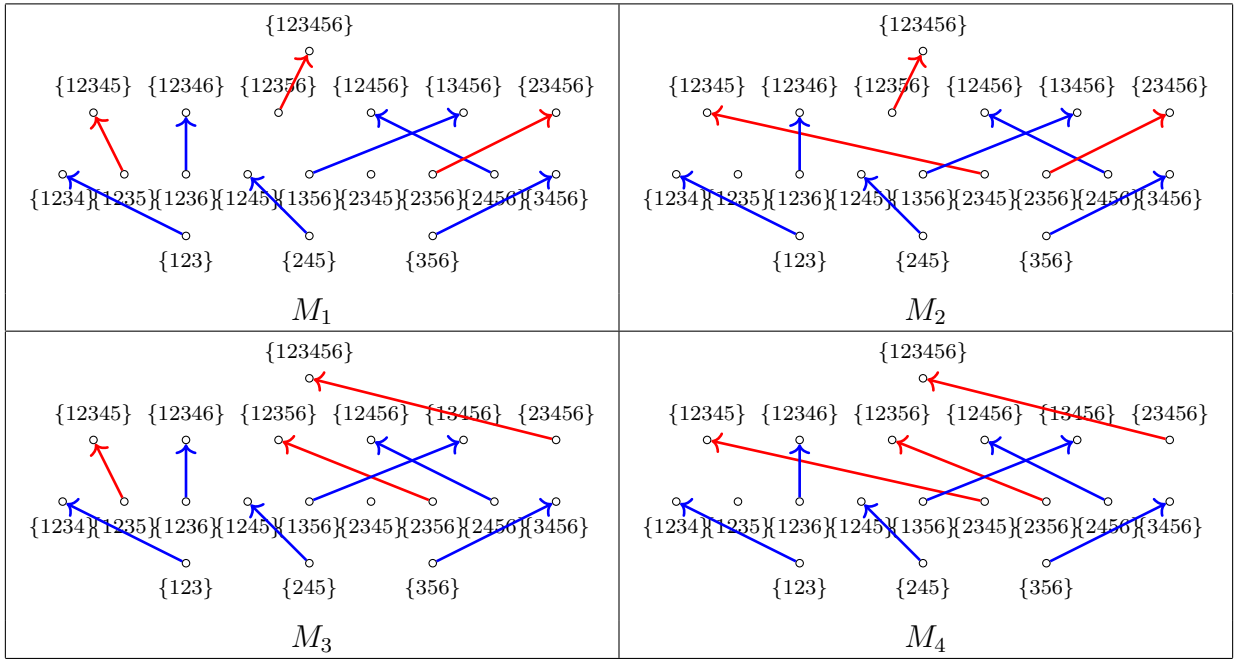

The homogeneous matchings $M_1,\ldots,M_4$ are shown in \Cref{fig:matchings}. Each of these matchings is, moreover, acyclic, because by~\cite[Lemma~3.3]{Roya}, a directed cycle will have to go between at least $6$ vertices all sharing the same monomial label.  
Note that for every $i \in \{1,\ldots,4\}$, the number of $M_i$-critical $j$-cells of $X$ are equal to the $j$-th Betti number of $I$. Thus, $X_{M_i}$ supports the minimal free resolution of $I$.

We now proceed to show that neither of the $X_{M_i}$ is polyhedral. Note that the two Morse complex $X_{M_1}$ and $X_{M_3}$ contain a single non-Scarf $3$-cell $\{2,3,4,5\}$, while $X_{M_2}$ and $X_{M_4}$ contain a single non-Scarf $3$-cell $\{1,2,3,5\}$. 

\begin{figure}[h]
\begin{tabular}{cc}
\begin{tikzpicture}[scale=.82]
\node [draw, circle, fill=white, inner sep=1pt, label=above:{\tiny{{${\{1245\}}$}}}] (1245) at (2.5,3) {};
\node [draw, circle, fill=white, inner sep=1pt, label=above:{\tiny{{${\{2345}\}$}}}] (2345) at (5.5,3) {};
\node [draw, circle, fill=black, inner sep=1pt, label=below:{\tiny{${\{124\}}$}}] (124) at (1,2) {};
\node [draw, circle, fill=black, inner sep=1pt, label=below:{\tiny{${\{125\}}$}}] (125) at (2,2) {};
\node [draw, circle, fill=black, inner sep=1pt, label=below:{\tiny{${\{145\}}$}}] (145) at (3,2) {};
\node [draw, circle, fill=white, inner sep=1pt,label=below:{\tiny{{${\{245\}}$}}}] (245) at (4,2) {};
\node [draw, circle, fill=black, inner sep=1pt, label=below:{\tiny{${\{234\}}$}}] (234) at (5,2) {};
\node [draw, circle, fill=black, inner sep=1pt, label=below:{\tiny{${\{235\}}$}}] (235) at (6,2) {};
\node [draw, circle, fill=black, inner sep=1pt,label=below:{\tiny{${\{345\}}$}}] (345) at (7,2) {};
\draw [line width=.5pt, ->] (1245)--(124);
\draw [line width=.5pt, ->] (1245)--(125);
\draw [line width=.5pt, ->] (1245)--(145);
\draw [color=blue,line width=1pt, <-] (1245)--(245);
\draw [line width=.5pt, ->] (2345)--(234);
\draw [line width=.5pt, ->] (2345)--(235);
\draw [line width=.5pt, ->] (2345)--(245);
\draw [line width=.5pt, ->] (2345)--(345);
\end{tikzpicture}
&
\begin{tikzpicture}[scale=.82]
\node [draw, circle, fill=white, inner sep=1pt, label=above:{\tiny{{${\{1234\}}$}}}] (1234) at (2.5,3) {};
\node [draw, circle, fill=white, inner sep=1pt, label=above:{\tiny{{${\{1235}\}$}}}] (1235) at (5.5,3) {};
\node [draw, circle, fill=black, inner sep=1pt, label=below:{\tiny{${\{234\}}$}}] (234) at (1,2) {};
\node [draw, circle, fill=black, inner sep=1pt, label=below:{\tiny{${\{134\}}$}}] (134) at (2,2) {};
\node [draw, circle, fill=black, inner sep=1pt, label=below:{\tiny{${\{124\}}$}}] (124) at (3,2) {};
\node [draw, circle, fill=white, inner sep=1pt,label=below:{\tiny{{${\{123\}}$}}}] (123) at (4,2) {};
\node [draw, circle, fill=black, inner sep=1pt, label=below:{\tiny{${\{125\}}$}}] (125) at (5,2) {};
\node [draw, circle, fill=black, inner sep=1pt, label=below:{\tiny{${\{135\}}$}}] (135) at (6,2) {};
\node [draw, circle, fill=black, inner sep=1pt,label=below:{\tiny{${\{235\}}$}}] (235) at (7,2) {};
\draw [line width=.5pt, ->] (1234)--(234);
\draw [line width=.5pt, ->] (1234)--(134);
\draw [line width=.5pt, ->] (1234)--(124);
\draw [color=blue,line width=1pt, <-] (1234)--(123);
\draw [line width=.5pt, ->] (1235)--(123);
\draw [line width=.5pt, ->] (1235)--(125);
\draw [line width=.5pt, ->] (1235)--(135);
\draw [line width=.5pt, ->] (1235)--(235);
\end{tikzpicture}\\
&\\
The $3$-cells of $\{2,3,4,5\}_{M_i}$ &
The $3$-cells of $\{1,2,3,5\}_{M_i}$\\
  $i=1,3$ &  $i=2,4$
\end{tabular}
\caption{Gradient paths}\label{fig:gradient}
\end{figure}
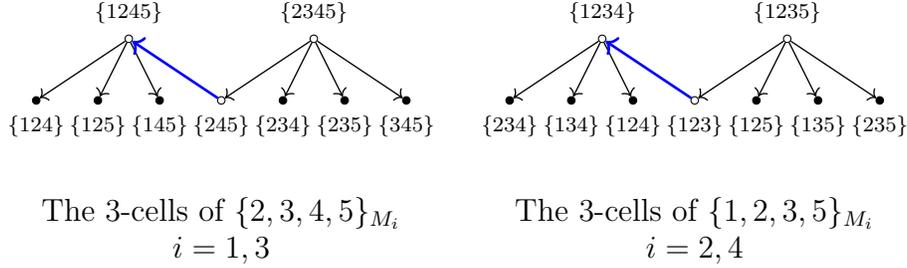

By \cref{lem:2-step-path} and following the gradient paths in \Cref{fig:gradient}, one can see the following. 

\begin{itemize}
\item  When $i=2,4$  and $\dim(\sigma)=2$,  if $\tau\neq \{1,2,3,5\}$, then $\sigma_{M_i}\subset \tau_{M_i}$ if and only if $\sigma\subset \tau$ by \Cref{lem:2-step-path}. 
However, when $\tau =\{1,2,3,5\}$,  we have  
$$
  \sigma_{M_i}\prec \{1,2,3,5\}_{M_i} 
  \iff 
  \sigma\in \big \{ 
   \{2,3,4\}, \{1,3,4\}, \{1,2,4\},\{1,2,5\},\{1,3,5\}, \{2,3,5\}
    \big \}
  $$ because of the gradient paths in \Cref{fig:gradient}.
Thus, 
\[
  \{1,5\}_{M_i}, \{1, 4\}_{M_i} \prec \{1,2,3,5\}_{M_i}\cap \{1,4,5,6\}_{M_i}
\]
and as these are maximal proper faces of  $\{1,2,3,5\}_{M_i}\cap \{1,4,5,6\}_{M_i}$ (note that $\{1, 4, 5\}_{M_i} \not \prec \{1,2,3,5\}_{M_i}$), we conclude  $\{1,2,3,5\}_{M_i}\cap \{1,4,5,6\}_{M_i}$ is not a face in $X_{M_i}.$ Hence $X_{M_i}$ is not a polyhedral cell complex.
To see why this gluing is not polyhedral, see also the bold dotted line in \Cref{fig:sub:nonpoly}.

\item When $i=1,3$ and $\dim(\sigma)=2$, if $\tau\neq \{2,3,4,5\}$, then $\sigma_{M_i}\prec \tau_{M_i}$ if and only if $\sigma\subset \tau$ by \cref{lem:2-step-path}. However, when $\tau =\{2,3,4,5\}$, 
  we have 
  $$
  \sigma_{M_i}\prec \{2,3,4,5\}_{M_i} 
  \iff 
  \sigma\in \big \{
  \{1,2, 4\},\{1,2, 5\},\{1,4,5\}, \{2,3,4\}, \{2,3,5\}, \{3,4,5\}
  \big \}
  $$ because of the gradient paths in \Cref{fig:gradient}.
Thus, 
\[
  \{3,4\}_{M_i}, \{1, 4\}_{M_i} \prec \{2,3,4,5\}_{M_i}\cap \{1,3,4,6\}_{M_i} 
\]
and as these are maximal proper faces of   $\{2,3,4,5\}_{M_i}\cap \{1,3,4,6\}_{M_i}$ (note that $\{1, 3, 4\}_{M_i} \not\prec \{2,3,4,5\}_{M_i}$), we conclude  $\{2,3,4,5\}_{M_i}\cap \{1,3,4,6\}_{M_i}$ is not a face in $X_{M_i}.$ Hence $X_{M_i}$ is not a polyhedral cell complex.   
\end{itemize}
\end{proof}

\begin{remarkbox}[{\bf An ideal with no minimal polyhedral resolution}]\label{r:minimal} 
    Note that more generally the ideal in \cref{nonexample} does not have a minimal polyhedral resolution. To see this, following the discussion earlier in this section, suppose $\omega$ is a $3$-polytope  with monomial label $x_1 \cdots x_{11}$ attached to $\Scarf(I)$
    along faces of $\Scarf(I)$, and the resulting complex $\MP$ is a polyhedral complex supporting a resolution of $I$. 
    Since $\MP$ consists of a single $3$-cell attached to a simplicial complex, the boundary of $\omega$ must be contained in $\Scarf(I)$ and hence simplicial; in other words  $\omega$ is a \say{simplicial polytope}. Then by Steinitz's lemma~\cite[Lemma 1.1]{Z2007}, if $f_i$ denotes the number of $i$-cells of $\omega$, we must have  
    \begin{equation}\label{e:Steinitz}
    f_1=3f_0-6 \qand f_2=2f_0-4. 
    \end{equation}
    Because of the monomial label $x_1 \cdots x_{11}$, the set of vertices of $\omega$ must contain $\{2,3,5\}$, and must be contained in $\{1,2,3,4,5\}$. Since $\omega$ is a $3$-cell, it has at least $4$ vertices, so following \eqref{e:the-taus}, $\omega$ falls into one of the three categories below. 
    \begin{enumerate}
        \item $\omega$ has vertex set $\{1,2,3,5\}$, in which case by \eqref{e:Steinitz} it has to be a simplex, but $\Scarf(I)$, whose facets are listed  in \eqref{e:Scarf-facets}, does not include the face $\{1,2,3\}$; a contradiction.
        \item $\omega$ has vertex set $\{2,3,4,5\}$, in which case it has to be a simplex again, but $\Scarf(I)$  does not contain the face $\{2,4,5\}$; a contradiction.
        \item $\omega$ has vertex set $\{1,2,3,4,5\}$, and so $f_0=5$, and by \eqref{e:Steinitz} $\omega$ must contain $f_1=9$ edges and $f_2=6$ triangles (simplicial $2$-cells).
             On the other hand the intersection of $\omega$ with every face of $\Scarf(I)$ must be a face.  This in particular means that if $\omega$ contains both vertices of an edge of $\Scarf(I)$, then it must contain the whole edge. It follows now from \eqref{e:Scarf-facets} that all $10$ edges between the vertices $1,\ldots,5$ of $\omega$ must be in $\omega$, which contradicts the requirements that $f_1=9$.
    \end{enumerate}
    So $\MP$ cannot be a polyhedral cell complex.

    Alternatively, one may use the fact that $\Scarf(I)$ is a simplicial complex, together with the fact that the new cell $\omega$ must intersect every face of $\Scarf(I)$ in a face of both $\Scarf(I)$ and $\omega$ to conclude that the $f$-vector of the boundary complex of $\omega$ is $(1,5,10,8,1)$. Computing the $h$-vector we get negative entries, which implies $(1,5,10,8,1)$ is not the $f$-vector of a simplicial polytope~\cite{S1980}, and we are done.
\end{remarkbox}

%%%%%%%%%%%%%%%%%%%%%%%%%%%%%%%%%%%%%%%%%%%%%%%%
%\bibliography{./ref.bib}{}
%\bibliographystyle{alpha}   

\newcommand{\etalchar}[1]{$^{#1}$}

%%%%%%%%%%%%%%%%%%%%%%%%%%%%%%%%%%%%%%%%%%%%%%%%
\end{document}